\documentclass[11pt]{article}

\usepackage[latin1]{inputenc}
\usepackage[english]{babel}
\usepackage{amsfonts,amssymb,amsmath}
\usepackage{mathrsfs}
\usepackage{mathtext}
\usepackage{textcomp}
\usepackage{graphicx}
\usepackage{epsfig}
\usepackage{rotating}

\usepackage{tikz}
\usetikzlibrary{arrows,decorations.pathmorphing,backgrounds,fit}

\usepackage[toc,page]{appendix}

\usepackage[dvips, colorlinks=true,pdfstartview=FitV, linkcolor=blue,citecolor=red,urlcolor=blue]{hyperref}

\newcommand{\lra}{\longrightarrow}
\newcommand{\sm}{\setminus}

\newcommand{\C}{\mathbb{C}}
\newcommand{\R}{\mathbb{R}}
\newcommand{\Z}{\mathbb{Z}}
\newcommand{\Q}{\mathbb{Q}}
\newcommand{\N}{\mathbb{N}}

\renewcommand{\S}{\mathbb{S}}

\renewcommand{\H}{\mathcal{H}}

\renewcommand{\P}{\mathrm{P}}

\newcommand{\Sp}{\mathbf{Sp}}
\newcommand{\SSp}{\mathbf{SSp}}
\newcommand{\SL}{\mathrm{SL}}
\newcommand{\Orb}{\mathcal{O}}

\newcommand{\Sig}{\Sigma}

\newcommand{\Id}{\mathrm{Id}}
\newcommand{\inter}{\mathrm{int}}

\newcommand{\D}{\mathrm{D}}
\newcommand{\Aa}{\mathbf{Area}}
\newcommand{\SpEmt}{(T_1,T_2,T_3,v_1,v_2)}
\newcommand{\OrbClo}{\overline{\SL(2,\R)\cdot\Psi(X_0)}}

\newcommand{\Hg}{\mathcal{H}(k_1,\dots,k_n)}
\newcommand{\Hgi}{\mathcal{H}_1(k_1,\dots,k_n)}

\newcommand{\dem}{\noindent {\textbf{Proof:} }}
\newcommand{\rem}{ \noindent \textbf{Remark: }}

\newcommand{\ie}{\textit{i.e. }}
\renewcommand{\geq}{\geqslant}
\renewcommand{\leq}{\leqslant}

\newcommand{\vide}{\varnothing}
\newcommand{\carre}{\hfill $\Box$\\}

\newcommand{\DS}{\displaystyle}

\newtheorem{theorem}{Theorem}[section]
\newtheorem{proposition}[theorem]{Proposition}
\newtheorem{prop-def}[theorem]{Proposition-definition}
\newtheorem{corollary}[theorem]{Corollary}
\newtheorem{lemma}[theorem]{Lemma}

\begin{document}

\title{\bf Parallelogram Decompositions and Generic Surfaces in $\H^\mathrm{hyp}(4)$}

\author{ DUC-MANH NGUYEN\\ IMB Bordeaux-Université Bordeaux 1\\ 351, Cours de la Libération\\ 33405 Talence Cedex\\FRANCE\\ \texttt{email: duc-manh.nguyen@math.u-bordeaux1.fr }}


\date{}

\maketitle

\begin{abstract}

The space $\H^\mathrm{hyp}(4)$ consists of pairs $(M,\omega)$, where $M$ is a hyper-elliptic Riemann surface of genus $3$, and $\omega$ is a holomorphic $1$-form having only one zero, which is located at a Weierstrass point of $M$. In this paper, we first show that {\em every} surface in $\H^\mathrm{hyp}(4)$ admits a decomposition into parallelograms and simple cylinders following a unique model. We then show that if this decomposition satisfies some specific condition, then the $\mathrm{GL}(2,\R)$-orbit of the surface is dense in $\H^\mathrm{hyp}(4)$. Using this criterion, we prove that there are generic surfaces in $\H^\mathrm{hyp}(4)$ with coordinates in any quadratic field, and there are Thurston-Veech surfaces with trace field of degree three over $\Q$ which are generic.

\end{abstract}

\section{Introduction}

Translation surfaces are flat surfaces with conical singularities and trivial linear holonomy, that is the holonomy of any closed curve is a translation in $\R^2$. The space of translation surfaces together with an oriented parallel line field is identified with the space of holomorphic $1$-forms on Riemann surfaces, which is stratified by the orders of the zeros of the $1$-form. Fix $g\geq 2$, if $k_1,\dots,k_n$ are some positive integers such that $k_1+\dots+k_n=2g-2$, we denote by $\Hg$ the moduli space of holomorphic $1$-forms on Riemann surfaces of genus $g$ which have exactly $n$ zeros with orders $(k_1,\dots,k_n)$. By a result of Kontsevich-Zorich \cite{KonZo}, we know that $\Hg$ has at most $3$ connected components. We denote by $\Hgi$ the subset of $\Hg$ consisting of surfaces of unit area.\\

\noindent There exists an action of $\SL(2,\R)$ on the space $\Hg$ which leaves invariant the Lebesgue measure, and preserves the set $\Hgi$. It is now a classical fact, due to Masur and Veech, that the $\SL(2,\R)$ action is ergodic in each component of $\Hgi$, a surface whose $\SL(2,\R)$-orbit is dense in its component is called {\em generic}. The $\SL(2,\R)$-orbit of almost all surfaces in each component is dense, however, the problem of determining whether the orbit of a particular surface is dense in its component is wide open. We only have a complete classification (due to McMullen and Calta,  \cite{McM2}, \cite{Cal}) for the case of genus $2$, where we have two strata, $\H(2)$ and $\H(1,1)$, each of which has a single connected component. Recall that the Veech group of a translation surface is the stabilizer subgroup for the action of $\SL(2,\R)$. It is a well-known fact that the $\SL(2,\R)$-orbit of a surface is a closed subset in its stratum if and only if its Veech group is a lattice of $\SL(2,\R)$. It turns out from the work of McMullen that, for translation surfaces of genus two, if the Veech group contains a hyperbolic element, then the $\SL(2,\R)$-orbit cannot be dense in the corresponding stratum.\\

More recently, Hubert-Lanneau-Möller (\cite{HubLanMo1}, \cite{HubLanMo2}) give some results on generic surfaces in the hyper-elliptic locus $\mathcal{L}$ of $\H^\mathrm{odd}(2,2)$, which is one of the two components of $\H(2,2)$. They show that, in contrast with the case of genus $2$,  there are generic surfaces in $\mathcal{L}$, that is the $\SL(2,\R)$-orbit is dense in $\mathcal{L}$, whose Veech group contains hyperbolic elements. Note that $\mathcal{L}$ is a closed, $\SL(2,\R)$-invariant subset of $\H^\mathrm{odd}(2,2)$, therefore, the closure of any $\SL(2,\R)$-orbit in $\mathcal{L}$ cannot exceed $\mathcal{L}$. The Thurston-Veech construction (\cite{Thur}, \cite{Vee89}) provides us with translation surfaces which are stabilized by some hyperbolic elements in $\SL(2,\R)$, these hyperbolic elements arise as products of parabolic elements. Hubert-Lanneau-Möller also show that there are surfaces in  $\mathcal{L}$ obtained from the Thurston-Veech construction whose $\SL(2,\R)$-orbit is dense in $\mathcal{L}$.\\

The stratum $\H(4)$ is the space of holomorphic $1$-form on Riemann surfaces of genus $3$ which have only one zero (the order is necessarily $4$). We have $\dim_\C\H(4)=6$, and $\H(4)$ has two connected components $\H^\mathrm{hyp}(4)$ and $\H^\mathrm{odd}(4)$ (see \cite{KonZo}). In this paper, we will be focusing on the connected component $\H^\mathrm{hyp}(4)$ which consists of holomorphic $1$-forms defined on hyper-elliptic Riemann surfaces. Equivalently, we can consider $\H^\mathrm{hyp}(4)$ as the space of translation surfaces of genus $3$ having only one singularity, such that there exists an isometric involution which  has exactly $8$ fixed points, and acts by $-\Id$ on the homology.\\

Before stating the main results of this paper, let us recall some basic definitions. On a translation surface, a {\em saddle connection} is a geodesic segment whose endpoints are singularities of the surface, which may coincide. For surfaces in $\H^\mathrm{hyp}(4)$, a saddle connection is then a geodesic loop joining the unique singularity to itself. If $\gamma$ is a saddle connection, we denote its length by $|\gamma|$. We can also associate to $\gamma$ together with a choice of orientation a vector $V(\gamma)\in \R^2$, which is the integral of the holomorphic $1$-form defining the flat metric along $\gamma$. In fact, the integral gives us a complex number, we view it as a vector in $\R^2$ by the standard identification $\C=\R\oplus \imath\R$. \\

\noindent Given a translation surface $\Sig$, a {\em cylinder} in $\Sig$ is an open subset which is isometric to the quotient $\R\times ]0;h[/\Z$, where $\Z$ is the cyclic group generated by $(x,y) \mapsto (x+\ell,y)$, and maximal with respect to this property. We will call $h$ the {\em height}, and $\ell$ the {\em width} of $C$, the modulus of $C$ is defined to be the ratio $\DS{h/\ell}$. Note that none of the parameters $h,\ell,m$ are invariant under $\SL(2,\R)$. By definition, we have a map from $\R\times ]0;h[$ to $\Sig$, which is locally isometric, with image $C$. This map can be extended by continuity to a map from $\R\times [0;h]$ to $\Sig$. We call the images of $\R\times\{0\}$ and $\R\times\{h\}$ under this map the {\em boundary components} of $C$. Each boundary component of $C$ is a concatenation of saddle connections, and freely homotopic to the simple closed geodesics in $C$. Remark that the two boundary components of $C$ are, in general, not disjoint subsets of $\Sig$, they can even coincide. We call $C$ a  {\em simple cylinder} when each of its boundary components consists of only one saddle connection.\\

A direction $\theta$ in $\S^1$ is said to be {\em completely periodic} if $\Sig$ is the union of the closures of the cylinders in this direction, in other words, any trajectory of the flow in this direction is either a closed geodesic or a saddle connection. 


\begin{theorem}\label{ThA}
On every surface in $\H^\mathrm{hyp}(4)$, there always exist four pairs of homologous saddle connections $\delta_i^\pm, i=1,\dots,4,$ such that

\begin{itemize}
\item[$\bullet$] $\delta_1^\pm$ bound a simple cylinder.

\item[$\bullet$] For $i=1,2,3, \delta_i^\pm$ and $\delta_{i+1}^\pm$ bound a topological disk, which is isometric to parallelogram in $\R^2$,

\item[$\bullet$] $\delta_4^\pm$ bound a simple cylinder.
\end{itemize}

\noindent The configuration of $\delta_1^\pm,\dots,\delta_4^\pm$ is shown in Figure 1.

\begin{figure}[!h]\label{DecompConfig}

\begin{center}
\begin{tikzpicture}[scale=0.4]
\draw[thin] (-2,1) -- (-3,-4)  (-1,-5) -- (4,-5) -- (5,0) -- (8,2) -- (8,6)  (9,9) -- (6,7)  (5,4) -- (0,4) -- (0,0) -- (5,0) -- (5,4) -- (8,6);

\draw[thin] (0,0) -- (-1,-5);

\draw[red!60!white, thin] (0,0) -- (-2,1) (-1,-5) -- (-3,-4) (5,4) -- (6,7) (8,6) -- (9,9);

\draw (7.5, 8) node[above] {$\delta_1^+$} (6.5,5) node[below] {$\delta_1^-$} (6.5,1) node[below] {$\delta_1^+$};
\draw (8,4) node[right] {$\delta_2^+$} (5,2) node[right] {$\delta_2^-$} (0,2) node[right] {$\delta_2^+$};
\draw (2.5,4) node[above] {$\delta_3^-$} (2.5,0) node[above] {$\delta_3^+$} (2.5,-5) node[below] {$\delta_3^-$};
\draw (4.5,-2.5) node[right] {$\delta_4^-$} (-0.5,-2.5) node[right] {$\delta_4^+$} (-2.5,-1.5) node[left] {$\delta_4^-$};

\draw[red!60!white, thin] (-2,-4.5) +(0,-0.2) -- +(0,0.2) (-1,0.5) +(0,-0.2) -- +(0,0.2);
\draw[red!60!white, thin] (5.5,5.5) +(-0.2,0.1) -- +(0.2,0.1) +(-0.2,-0.1) -- +(0.2,-0.1) (8.5,7.5) +(-0.2,0.1) -- +(0.2,0.1) +(-0.2,-0.1) -- +(0.2,-0.1);

\foreach \x in {(0,0), (5,0), (5,4), (0,4), (8,2), (8,6) , (9,9), (6,7), (4,-5), (-1,-5), (-3,-4), (-2,1)}{
\filldraw[gray] \x circle (2pt);}

\end{tikzpicture}
\caption{Decomposition of surfaces in $\H^\mathrm{hyp}(4)$ into parallelograms and simple cylinders}
\end{center}

\end{figure}
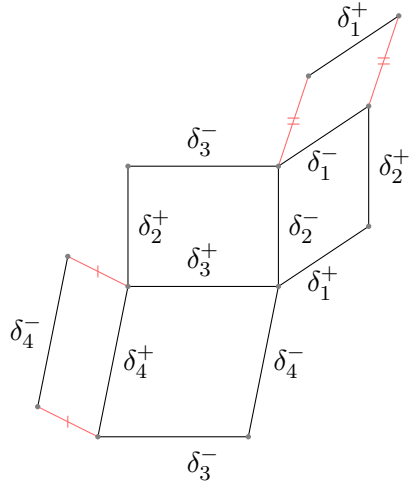

\end{theorem}

Let $\Sig_0$ be a surface in $\H^\mathrm{hyp}_1(4)$, and $\delta_i^\pm, i=1,\dots,4,$ be as in Theorem \ref{ThA}. Cutting $\Sig$ along $\delta_3^\pm$, we get two connected components whose boundary consists of two geodesic segments. Gluing those geodesic segments together, we then get a flat torus, which will be denoted by $\Sig'$, and a surface in $\H(2)$. On the torus $\Sig'$, we denote the geodesic segment corresponding to $\delta_3^\pm$ by $\delta_3$. As a subsurface of $\Sig$, $\Sig'$ inherits a parallel line field, therefore we can view it as a pair $(M,\omega)$, where $M$ is a Riemann surface of genus one, and $\omega$ is a non-zero holomorphic $1$-form on $M$. Equivalently, we can  identify $\Sig'$ with the quotient $\R^2/\Lambda$, where $\Lambda$ is a lattice in $\R^2$, which is the image of the map $\DS{H_1(M,\Z) \lra \C\simeq \R^2: c \mapsto \int_c\omega}$. A vector in $\R^2$ is said to be {\em generic} with respect to $\Lambda$ if it is not collinear with any vector in $\Lambda$. We have

\begin{theorem}\label{ThB}
Suppose that $\delta_1^\pm$ and $\delta_3^\pm$ are parallel, that is $V(\delta_1^\pm)$ and $V(\delta_3^\pm)$ are collinear, and $V(\delta_3)=V(\delta_3^\pm)$ is generic with respect to $\Lambda$, then $\SL(2,\R)\cdot\Sig_0$ is dense in $\H_1^\mathrm{hyp}(4)$.
\end{theorem}

Using this result, we obtain

\begin{corollary}\label{CorB}
Let $\Sig_0$ be a surface in $\H_1^\mathrm{hyp}(4)$. Suppose that the horizontal direction is completely periodic for $\Sig_0$, and that $\Sig_0$ is decomposed into three horizontal cylinders whose moduli are independent over $\Q$, then  $\SL(2,\R)\cdot\Sig_0$ is dense in $\H_1^\mathrm{hyp}(4)$.
\end{corollary}

The proof of Theorem \ref{ThA} relies on the action of the hyper-elliptic involution on the surfaces in $\H^\mathrm{hyp}(4)$. The key  ingredient of the proof is Lemma \ref{InvSC}, which says that, on a translation surface of genus one or two, any saddle connection invariant under the distinguished involution of the surface is contained in a simple cylinder.\\
To prove Theorem \ref{ThB}, we will show that the orbit closure contains all the surfaces admitting a splitting as in Theorem \ref{ThA} with $\delta_2^\pm$ parallel to $\delta_3^\pm$. Consequently, the orbit closure contains all the Veech surfaces, and in particular all the square-tiled surfaces. Since the set of square-tiled surfaces is dense in $\H_1^\mathrm{hyp}(4)$, we deduce that the orbit closure is the whole component. The proof of Theorem \ref{ThB} uses a theorem of Ratner on action of unipotent subgroups on homogeneous spaces.\\
To prove Corollary \ref{CorB}, we prove that one can find in the $\SL(2,\R)$-orbit closure of $\Sig_0$ a surface which satisfies the condition of Theorem \ref{ThB}. It is easy to construct surfaces in $\H^\mathrm{hyp}(4)$ with coordinates in a quadratic field over $\Q$ which satisfy the condition of Theorem \ref{ThB}, therefore, we have an affirmative answer to a question in \cite{HubLanMo3} (see Section \ref{ExSect}). We will also construct explicitly some Thurston-Veech surfaces with trace field of degree three over $\Q$ which satisfy the hypothesis of Corollary \ref{CorB}.\\

\section{Simple cylinder invariant under the involution}

\subsection{Translation surfaces of genus one}

Translation surfaces of genus one are simply flat tori. We denote by $\H(0)$ (resp. $\H(0,0)$) the space of triples $(M,\omega,P)$ (resp.
quadruplet $(M,\omega,P_1,P_2)$), where $M$ is a Riemann surface of genus one, $\omega$ is a nonzero holomorphic $1$-form on $M$, and $P$ (resp.
$P_1$ and $P_2$) is a marked point (resp. are marked points) of $M$. In both cases, we will call the lattice in $\R^2$ obtained by integrating $\omega$ along elements of $H_1(M,\Z)$ the {\it associated lattice} of the considered translation surface. If $\Sig$ is an element of $\H(0)$ or $\H(0,0)$, we denote by $\Lambda(\Sig)$ the lattice associated to $\Sig$.\\ 


\noindent Note that the holomorphic $1$-form determines a flat metric structure together with a choice of vertical direction at every point of the surface. For each surface in $\H(0)$, and $ \H(0,0)$, we have a distinguished isometric involution which acts like $-\Id$ on the homology of the surface, and either fixes the unique marked point (in the case of $\H(0)$), or exchanges the two marked points (in the case of $\H(0,0)$). As usual, we will call a geodesic segment joining marked points a {\em saddle connection}. In the case of $\H(0)$, a saddle connection is just a simple closed geodesic passing through the marked point.\\

\subsection{Saddle connection preserved by the involution}

Let us prove the following lemma, which is the key ingredient for the proof of Theorem \ref{ThA},

\begin{lemma}\label{InvSC}
Let $\gamma$ be a saddle connection on a translation surface $\Sig$ which belongs to one of the following strata $\H(0), \H(0,0), \H(2), \H(1,1)$. Suppose that $\gamma$ is invariant under the distinguished involution in the cases $\H(0)$ and $\H(0,0)$, or under the hyper-elliptic involution in the cases $\H(2)$ and $\H(1,1)$, then there exists a pair of saddle connections $(\eta^+,\eta^-)$ which bound a simple cylinder $C$ containing $\gamma$, \ie

\begin{itemize}

\item[.] $\overline{C}\sm C= \eta^+\cup \eta^-$,

\item[.] $\inter(\gamma) \subset C$.

\end{itemize}

In the case $\H(0)$, actually $\eta^+\equiv \eta^-$, in all others case $\eta^+$ and $\eta^-$ are distinct.
\end{lemma}

\dem We will prove this lemma case by case.\\

\noindent \underline{\bf Case $\H(0)$:} \\
\noindent In this case $\gamma$ is a simple closed geodesic passing through the marked point. Let $\eta$ be  any simple closed geodesic which meets  $\gamma$ only at the marked point, then we can take $\eta^+=\eta^-=\eta$.\\

\noindent \underline{\bf Case $\H(0,0)$:} \\
\noindent In this case $\gamma$ is a geodesic segment joining two marked points $P_1, P_2$ of $\Sig$. Using the action of $\SL(2,\R)$, we can assume that $\gamma$ is horizontal. Let $\Psi_t, t\in \R,$ denote the vertical flow on $\Sig$. There exists a minimal value $t_0>0$ such that $\Psi_{t_0}(\gamma)\cap \gamma \neq \vide$. Observe that $\Psi_{t_0}(\gamma)$ must contain one endpoint of $\gamma$, without loss of generality, we can assume that $P_1\in \Psi_{t_0}(\gamma)$.\\
\noindent  By the definition of $t_0$, we have an isometric immersion $\Phi$ from the rectangle $R=[0;|\gamma|]\times [0;t_0]$ into $\Sig$ whose restriction into $\inter(R)$ is an embedding. We can suppose that $\Phi$ maps the lower side of $R$ onto $\gamma$. Let $\tilde{P}_i^b, i=1,2,$ denote the two endpoints of the lower side of $R$ so that $\Phi(\tilde{P}_i^b)=P_i$. By assumption, there exists a point  $\tilde{P}_1^t$ in the upper side of $R$ such that $\Phi(\tilde{P}^t_1)=P_1$. Let $\tilde{\eta}$ denote the geodesic segment in $R$ joining $\tilde{P}_1^b$ to $\tilde{P}_1^t$, then $\eta^+=\Phi(\tilde{\eta})$ is a simple closed geodesic in $\Sig$ which meets $\gamma$ only at $P_1$. Let $\eta^-$ be the image of $\eta^+$ under the distinguished involution of $\Sig$, we see that $\eta^-$ is  parallel to $\eta^+$, and meets $\gamma$ only at $P_2$. It is easy to check that $\eta^+$ and $\eta^-$ cut $\Sig$ into two cylinders, one of which contains  $\gamma$.\\

\noindent \underline{\bf Case $\H(2)$:} \\
\noindent Let $P$ denote the unique singularity of $\Sig$, and $\tau$ denote the hyper-elliptic involution of $\Sig$. In this case $\gamma$ is a geodesic segment joining $P$ to itself, and invariant under $\tau$. Note that $\tau$ reverse the orientation of $\gamma$, and since $\tau(P)=P$, it also fixes the midpoint of $\gamma$.\\
\noindent We can assume that $\gamma$ is horizontal. As before, let $\Psi_t, t\in \R,$ denote the vertical flow on $\Sig$. The
same argument as in the previous case shows that we have an immersion $\Phi$ from a rectangle $R \subset \R^2 $ into $\Sig$ such that $\Phi_{|\inter(R)}$ is an embedding, $\Phi$ maps the lower side of $R$ onto $\gamma$, and there exists a point $\tilde{P}$ in the upper side of $R$ which is mapped to $P$. Let $\tilde{\Delta}$ denote the triangle whose vertices are $\tilde{P}$ and the two endpoints of the lower side of $R$. Since $\tilde{\Delta}$ is contained in $R$, the restriction $\Phi_{|\inter(\tilde{\Delta})}$ is an embedding, moreover the images of the sides of  $\tilde{\Delta}$ by $\Phi$ are three distinct saddle connections, which meet one another only at $P$. Therefore,  $\Delta=\Phi(\tilde{\Delta})$ is an embedded triangle in $\Sig$ whose vertices coincide with $P$. By construction, $\gamma$ is a side of $\Delta$, let $\sigma_1,\sigma_2$ denote the two other sides. Let $\Delta', \sigma'_1,\sigma'_2$ denote the images of $\Delta,\sigma_1,\sigma_2$ under $\tau$ respectively. Observe that $\Delta'$ is also an embedded triangle in $\Sig$, and $\gamma$ is a common side of $\Delta'$ and $\Delta$. Here we have two possibilities:

\begin{itemize}

\item[$\bullet$] $\Delta$ and $\Delta'$ have another common side other than $\gamma$, that is, either $\sigma'_1=\sigma_1$, or $\sigma_2=\sigma'_2$. In this case $\Delta\cup\Delta'$ is a simple cylinder, and we are done.

\item[$\bullet$] $\gamma$ is the only common side of $\Delta$ and $\Delta'$. In this case, $\Delta\cup\Delta'$ is an embedded parallelogram in $\Sig$. Let us show that $\sigma_1$ and $\sigma'_1$ bound a cylinder disjoint from $\Delta\cup\Delta'$.  Recall the the cone angle at $P$ is $6\pi$, and the action of $\tau$ at $P$ is the rotation of angle $3\pi$.  Fix an orientation for $\gamma$, consider $\gamma$ as a part of $\partial\Delta$ (resp. $\partial\Delta'$), we then have an orientation for $\sigma_1,\sigma_2$ (resp. $\sigma'_1,\sigma'_2$) subsequently. Consider a small disk $\D$ centered at $P$. The intersection  of any oriented saddle connection with $\D$ is the union of an outgoing ray, and an incoming ray. These two rays specify a pair of angles at $P$, since $\Sig$ is a translation surface, this pair of angles is either $(\pi, 5\pi)$, or $(3\pi,3\pi)$.  Since $\gamma$ is invariant under $\tau$, the pair of angle specified by $\gamma$ is $(3\pi,3\pi)$, meanwhile the pair of angles specified by $\sigma_1$ is $(\pi,5\pi)$ since $\sigma'_1=\tau(\sigma_1) \neq \sigma_1$. We claim that the outgoing and the incoming rays of $\sigma_1$ are contained in the same half disk cut out by the outgoing and the incoming rays of $\gamma$. Indeed, suppose that the outgoing and the incoming rays of $\sigma_1$ do not belong to the same half disk (see Figure 2, Case a)), then by considering the sum of the angles in $\Delta$, we see that the pair of angles specified by $\sigma_2$ is $(3\pi,3\pi)$, which means that $\sigma_2=\tau(\sigma_2)=\sigma'_2$, but this is excluded by the hypothesis.

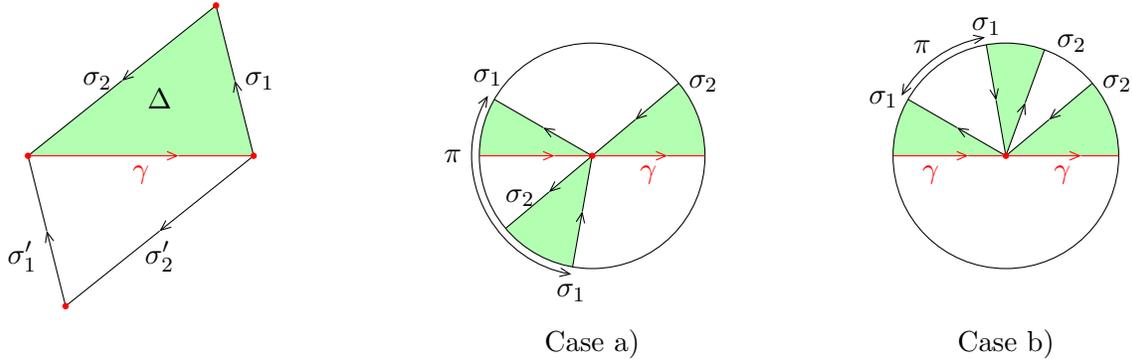
\begin{figure}[!ht]\label{Fig2}
\begin{center}
\begin{tikzpicture}[scale=0.5]
\fill[green!30!white] (-12,0) -- (-6,0) -- (-7,4) -- cycle;
\draw[red, thin, >=angle 45, ->] (-12,0) -- (-8,0); \draw[thin, red] (-8,0) -- (-6,0);
\draw[>= angle 45, ->] (-6,0) -- (-6.5,2); \draw[>=angle 45, ->] (-7,4) -- (-9.5,2); \draw[>=angle 45,->](-11,-4) -- (-11.5,-2);
\draw[>= angle 45, ->] (-6,0) -- (-8.5,-2);
\draw (-6.5,2) -- (-7,4) (-9.5,2) -- (-12,0) -- (-11.5,-2) (-11,-4) -- (-8.5,-2);
\filldraw[red] (-12,0) circle (2pt) (-6,0) circle (2pt) (-7,4) circle (2pt) (-11,-4) circle (2pt);
\draw[red] (-9,0) node[below] {$\gamma$}; \draw (-6.5,2) node[right] {$\sigma_1$} (-9.5,2) node[left] {$\sigma_2$};
\draw (-11.5,-2) node[below left] {$\sigma'_1$} (-8.5,-2) node[below] {$\sigma'_2$};
\draw (-8.5,1) node[above] {$\Delta$};

\fill[green!30!white] (3,0) -- +(150:3cm) arc  (150:180:3cm) -- cycle;
\fill[green!30!white] (3,0) -- +(0:3cm) arc (0:40:3cm) -- cycle;
\fill[green!30!white] (3,0) -- +(220:3cm) arc (220:260:3cm) -- cycle;
\draw (3,0) circle (3cm); 
\draw[red, >= angle 45, ->] (0,0) -- (2,0); \draw[red, >=angle 45, ->] (3,0) --(5,0);
\draw[>= angle 45, ->] (3,0) -- +(150:1.5cm); \draw[>= angle 45, ->] (3,0) +(40:3cm) -- +(40: 1.5cm); \draw[>= angle 45, ->] (3,0) +(260:3cm) -- +(260: 1.5cm); \draw[>= angle 45, ->] (3,0) -- +(220:1.5cm);
\draw[red] (2,0) -- (3,0) (5,0) -- (6,0);
\draw (3,0) +(150:1.5cm) -- +(150:3cm) +(0,0) -- +(40:1.5cm) +(0,0) -- +(260:1.5cm) +(220:1.5cm) -- +(220:3cm);
\filldraw[red] (3,0) circle (2pt);
\draw[>=angle 45, <->] (3,0)  +(150:3.2cm) arc (150: 260:3.2cm);
\draw[red] (4.5,0) node[below] {$\gamma$};
\draw (3,0) +(40:3cm) node[right] {$\sigma_2$} +(150:3.2cm) node[above] {$\sigma_1$} +(260:3.2cm) node[below] {$\sigma_1$} +(220:2.5cm) node[above] {$\sigma_2$} +(180:3.2cm) node[left] {$\pi$};

\fill[green!30!white] (14,0) -- +(150:3cm) arc  (150:180:3cm) -- cycle;
\fill[green!30!white] (14,0) -- +(0:3cm) arc (0:40:3cm) -- cycle;
\fill[green!30!white] (14,0) -- +(70:3cm) arc (70:100:3cm) -- cycle;
\draw (14,0) circle (3cm); 
\draw[red, >= angle 45, ->] (11,0) -- (13,0); \draw[red, >=angle 45, ->] (14,0) --(16,0);
\draw[>= angle 45, ->] (14,0) -- +(150:1.5cm); \draw[>= angle 45, ->] (14,0) +(40:3cm) -- +(40: 1.5cm); \draw[>= angle 45, ->] (14,0) +(100:3cm) -- +(100: 1.5cm); \draw[>= angle 45, ->] (14,0) -- +(70:1.5cm);
\draw[red] (13,0) -- (14,0) (16,0) -- (17,0);
\draw (14,0) +(150:1.5cm) -- +(150:3cm) +(0,0) -- +(40:1.5cm) +(0,0) -- +(100:1.5cm) +(70:1.5cm) -- +(70:3cm);
\filldraw[red] (14,0) circle (2pt);
\draw[>=angle 45, <->] (14,0) +(100:3.2cm) arc (100:150:3.2cm);
\draw[red] (12,0) node[below] {$\gamma$} (15.5,0) node[below] {$\gamma$};
\draw (14,0) +(40:3cm) node[right] {$\sigma_2$} +(70:3.2cm) node[ right] {$\sigma_2$} +(100:3cm) node[above] {$\sigma_1$} +(150:3cm) node[left] {$\sigma_1$} +(125:3cm) node[above left] {$\pi$};

\draw (3,-5) node {Case a)} (14,-5) node {Case b)};
\end{tikzpicture}

\caption{Configurations of geodesics rays at $P$}
\end{center}

\end{figure}

We know that the action of $\tau$ on $H_1(\Sig,\Z)$ is $-\Id$, which implies that $ \sigma_1-\sigma'_1=0 \text{ in } H_1(\Sig,\Z)$. It follows that $\sigma_1$ and $\sigma'_1$ cut $\Sig$ into two connected components, each of which is equipped with a flat metric structure with piecewise geodesic boundary. Consider the connected component which does not contain $\gamma$. This component does not contain any singularity in its
interior, and since the angle between the two rays of $\sigma_1$ at $P$ measured inside this component is $\pi$, we deduce that there is no
singularities in its boundary. The only flat surface with two geodesic boundary components with no singularities is a cylinder. Therefore, we can conclude that $\sigma_1$ and $\sigma'_1$ bound a cylinder $C$ disjoint from $\Delta\cup\Delta'$.\\
Consider the subsurface $\Sig'=\Delta\cup\Delta'\cup C$ of $\Sig$. We first observe that $\Sig'$ is invariant under $\tau$. Topologically,
$\Sig'$ is the complement in a torus of two open disks whose boundaries meet at one point. We can construct $\Sig'$ by gluing two parallelograms
so that the restriction of $\tau$ into $\Sig'$ is realized by the central symmetries in both parallelograms. Elementary geometry shows that one can find a saddle connection $\eta^+$ in $\Delta\cup C$ which crosses $\sigma_1$ once. Let $\eta^-$ denote the image of $d$ under $\tau$, then $\eta^+$ and $\eta^-$ bound  a simple cylinder containing $\gamma$. Since $\tau$ preserves $\gamma$ and reverses its orientation, we see that $\tau$ preserves the cylinder bounded by $\eta^+$ and $\eta^-$, and $\tau(\eta^+)=\eta^-, \tau(\eta^-)=\eta^+$.

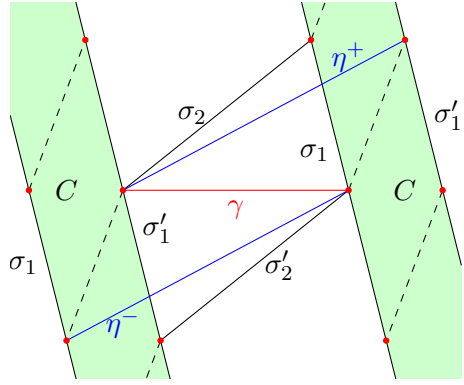
\begin{figure}[!ht]\label{Fig3}
\begin{center}
\begin{tikzpicture}[scale=0.5]
\clip (-6,-5) rectangle (6,5);
\fill[green!20!white] (4.5,-6) -- (7,-6) -- (4,6) -- (1.5,6) -- cycle;
\fill[green!20!white] (-4.5,6) -- (-7,6) -- (-4,-6) -- (-1.5,-6) -- cycle;
\draw[red] (-3,0) -- (3,0);
\draw (2,4) -- (-3,0)  (-2,-4) -- (3,0) (6,-2) -- (4,6) (-6,2) -- (-4,-6);
\draw (-4.5,6) -- (-1.5,-6) (1.5,6) -- (4.5,-6);
\draw[blue] (-3,0) -- (4.5,4) (3,0) -- (-4.5,-4);
\draw[dashed] (-3,0) -- (-4.5,-4) (3,0) -- (4.5,4) (-5.5,0) -- (-4,4) (5.5,0) -- (4,-4) (-2,-4) -- (-3.5,-8) (2,4) -- (3.5,8);

\filldraw[red] (-3,0) circle (2pt) (3,0) circle (2pt) (2,4) circle (2pt) (-2,-4) circle (2pt) (-4.5,-4) circle (2pt) (-5.5,0) circle (2pt) (4.5,4) circle (2pt) (5.5,0) circle (2pt) (-4,4) circle (2pt) (4,-4) circle (2pt);

\draw[red] (0,0) node[below] {$\gamma$};
\draw (-0.5,2) node[left] {$\sigma_2$} (0.5,-2) node[right] {$\sigma'_2$} (2.75,1) node[left] {$\sigma_1$} (-2.75,-1) node[right] {$\sigma'_1$} (5,2) node[right] {$\sigma'_1$} (-5,-2) node[left] {$\sigma_1$};

\draw[blue] (3,4.2) node[below] {$\eta^+$} (-3, -4.2) node[above] {$\eta^-$};

\draw (-4.5,0) node {$C$} (4.5,0) node {$C$};

\end{tikzpicture}
\caption{Existence of $\eta^\pm$}
\end{center}
\end{figure}

\end{itemize}

\noindent \underline{\bf Case $\H(1,1)$:}\\
Let $\{P_1,P_2\}$ denote the singularities of $\Sig$, the cone angles at both $P_1$ and $P_2$ are $4\pi$. Recall that in this case, the hyper-elliptic involution $\tau$ exchanges $P_1$ and $P_2$, therefore $\gamma$ must be a saddle connection joining $P_1$ to $P_2$. Without loss of generality, we can assume that $\gamma$ is horizontal.  As we have seen in the previous cases, there exists an embedded triangle $\Delta$ in $\Sig$ bounded by $\gamma$ and two other saddle connections $\sigma_1$ and $\sigma_2$. Since there are only two singularities, one of the two sides $\sigma_1$ and $\sigma_2$ must joint a singularity to itself, therefore we can assume that $\sigma_1$ joins $P_1$ to itself.\\
Let $\Delta', \sigma'_1,\sigma'_2$ denote the images of $\Delta,\sigma_1,\sigma_2$ under $\tau$ respectively. Since $\tau$ exchanges $P_1$ and $P_2$, $\sigma_1$ and $\sigma'_1$ are two distinct saddle connections. We choose the orientation for $\gamma$ to be from $P_1$ to $P_2$, and choose the orientation of $\sigma_1$ and $\sigma_2$ (resp. $\sigma'_1$ and $\sigma'_2$) coherently to get an orientation for  the boundary of $\Delta$ (resp. $\Delta'$). Consider two small disks $\D_1,\D_2$ centered at $P_1,P_2$ respectively. The intersection of $\sigma_1$ with $\D_1$ consists of an outgoing ray and an incoming ray, while the intersection of $\gamma$ with $\D_1$ consists of only an out going ray. Let $\theta$ be the angle between the outgoing and the incoming rays of $\sigma_1$ measured along the sector of $\D_1$ that does not contain $\gamma$. We have two cases:

\begin{itemize}
\item[.] $\theta=3\pi:$ In this case, the angle between the two rays of $\sigma_1$ measured along the other sector of $\D_1$ is $\pi$. A simple
computation of angles shows that we must have $\sigma_2=\sigma'_2$ as subset of $\Sig$, which implies that $\Delta\cup\Delta'$ is actually a cylinder invariant under $\tau$ and bounded by $\sigma_1$ and $\sigma'_1$, and the lemma follows immediately.

\item[.] $\theta=\pi:$ Since $\sigma_1-\sigma'_1=0$ in $H_1(\Sig,\Z)$, by cutting $\Sig$ along $\sigma_1$ and $\sigma'_1$, we obtain two flat surfaces with piecewise geodesic boundary. Observe that the component which does not contain $\gamma$ has no singularities in the interior, and
since the angle between the two rays of $\sigma_1$ measured inside this surface is $\pi$, we see that it has no singularities in the boundary. It
follows that this component is a cylinder $C$ bounded by $\sigma_1$ and $\sigma'_1$. Now, using the same argument as in the case $\H(2)$, we see that there exists a pair of saddle connections $\eta^\pm$ in $\Delta\cup\Delta'\cup C$ which are exchanged by $\tau$, and bound a simple cylinder containing $\gamma$.
\end{itemize}

\carre

\section{Proof of Theorem \ref{ThA}}

\subsection{Existence of simple cylinder on hyper-elliptic translation surfaces}

To prove Theorem \ref{ThA}, we first show

\begin{lemma}\label{SCylEx}

For any $g\geq 2$, on every surface of the stratum $\H^{\mathrm{hyp}}(2g-2)$, there always exists a simple cylinder which is invariant under the
hyper-elliptic involution.

\end{lemma}

\dem Let $\Sig$ be a surface in the stratum $\H^\mathrm{hyp}(2g-2)$. A construction due to Veech (see \cite{Vee89}, \cite{HubLanMo3}) allows us to
construct $\Sig$  from a $2g$-gon $\P$ in $\R^2$ centered at the origin, and invariant under the central symmetry of $\R^2$. The polygon $\P$ is
not necessarily convex, however it  has a horizontal diagonal $d$ which passes through the origin contained in the interior. Let $A_0,B_0$
denote the left and right endpoints of $d$ respectively. We denote by $A_1,\dots,A_{2g-1}$ (resp. $B_1,\dots,B_{2g-1}$) the vertices of $\P$
above (resp. below) the diagonal $d$ in the counter-clockwise order. We consider by convention that $A_{2g}=B_0$, and $B_{2g}=A_0$. The surface
$\Sig$ is obtained by identifying the opposite sides of $\P$.\\
Let $y:\R^2\lra \R$ denote the vertical coordinate function of $\R^2$. Let $i_0$ be the smallest index in $\{0,1,\dots, 2g-1\}$ so that
$y(A_{i_0})= \max\{y(A_0), \dots, y(A_{2g-1})\}$. Note that we have $0<i_0<2g$ since $y(A_0)=y(B_0)=0$. By the choice of $i_0$, we see that the diagonal $\overline{A_{i_0-1}A_{i_0+1}}$ is contained inside $\P$. By symmetry, the diagonal $\overline{B_{i_0-1}B_{i_0+1}}$ is also contained inside $\P$. Since the sides $\overline{A_{i_0-1}A_{i_0}}$ and $\overline{A_{i_0}A_{i_0+1}}$ are identified with $\overline{B_{i_0-1}B_{i_0}}$ and
$\overline{B_{i_0}B_{i_0+1}}$ respectively, it follows that the union of the two triangles $\Delta^u=(A_{i_0-1}A_{i_0}A_{i_0+1})$ and
$\Delta^l=(B_{i_0-1}B_{i_0}B_{i_0+1})$ is projected to a simple cylinder $C$ of $\Sig$. Now,  the hyper-elliptic involution of $\Sig$
corresponds to the central  symmetry at the origin, which interchanges the two triangles $\Delta^u$, and $\Delta^l$, therefore the
hyper-elliptic involution preserves $C$, and exchanges its two boundary components .\carre

\rem This lemma is also true for surfaces in $\H^{\mathrm{hyp}}(g-1,g-1)$.

\subsection{Proof of Theorem \ref{ThA}}

Let $\Sig$ be a surface in $\H^\mathrm{hyp}(4)$, we denote by $\tau$ the hyper-elliptic involution of $\Sig$. By Lemma \ref{SCylEx}, we know
that there exists a simple cylinder $C_1$ in $\Sig$ bounded by a pair of saddle connections $(\delta_1^+,\delta_1^-)$ such that $\tau(C_1)=C_1$ and $\tau(\delta_1^+)=\delta_1^-$. Cutting off $C_1$ from $\Sig$, we then get a surface whose boundary is an eight figure, \ie the union of two circles meeting at one point. Splitting the common point of the two circles into two points gives us two geodesic segments (corresponding to the pair
$(\delta_1^+,\delta_1^-)$), gluing these two segments together, we then get a surface $\Sig'$ in $\H(1,1)$ with a marked saddle connection which
will be denoted by $\delta_1$.\\
Since $\tau$ preserves $C_1$, and exchanges $\delta^+_1$ and $\delta^-_1$, its restriction $\tau'$ to $\Sig'$ is the hyper-elliptic
involution of $\Sig'$, and preserves the saddle connection $\delta_1$. By Lemma \ref{InvSC}, we know that there exists a pair of saddle
connections $(\delta_2^+,\delta_2^-)$ in $\Sig'$ which bound a simple cylinder $C_2$ containing $\delta_1$. Again,  we have that $\tau'$
preserves $C_2$ and exchanges $\delta_2^+$ and $\delta_2^-$. Note that since $\delta_2^+$ and $\delta_2^-$ meet $\delta_1$ at only the endpoints
of $\delta_1$, which are the singularities of $\Sig'$, we deduce that $\delta_2^+$ and $\delta_2^-$ are a pair of homologous saddle
connections in the initial surface $\Sig$.\\
Now, cut off $C_2$ from $\Sig'$, what is left is a surface with two boundary components corresponding to $\delta_2^+$ and $\delta_2^-$.
Gluing the two boundary components so that the two singularities are identified, we get a surface in $\H(2)$ with a marked saddle connection,
which is invariant by the hyper-elliptic involution. Lemma \ref{InvSC} then allows us to continue the procedure until we are left with a simple
cylinder. Since in each step, we cut out a simple cylinder, a simple computation on Euler character shows that we get to this situation after
four steps. The result of this procedure is that we have found four pairs of homologous saddle connections $(\delta_i^+,\delta_i^-), i=1,\dots,4,$
in $\Sig$ which satisfy the properties asserted in the statement of the theorem. \carre

\begin{corollary}\label{SepSC}

There exists on any surface $\Sig$ in $\H^\mathrm{hyp}(4)$ a pair of homologous saddle connections which are exchanged by the hyper-elliptic involution, and  decompose $\Sig$ into a union of a surface in $\H(2)$, and a surface in $\H(0,0)$. In both components of this decomposition, this pair of saddle connections corresponds to a saddle connection invariant under the (distinguished) involution.

\end{corollary}

\dem Let $(\delta_i^+,\delta_i^-), i=1,\dots,4,$ be the saddle connections in $\Sig$ satisfying the properties in Theorem \ref{ThA}. It is easy
to check that both pairs $(\delta^+_2,\delta_2^-)$ and $(\delta_3^+,\delta_3^-)$ satisfy the property asserted in the corollary.\carre

%
%
%
%
%
%
%
%

\section{Splitting of surfaces in $\H^\mathrm{hyp}(4)$}

\subsection{Flat torus with a marked geodesic segment}

Throughout this paper, by a 'flat torus' we will mean a Riemann surface of genus one together with a non-zero holomorphic $1$-form.
Equivalently, we identify a flat torus with the quotient $\C/\Lambda$, where $\Lambda$ is a lattice isomorphic to $\Z\oplus\Z$. Using this identification, we can associate to any oriented geodesic $s$ in the torus a vector $V(s)\in\R^2$. If $u=(x_1,y_1)$ and $v=(x_2,y_2)$ are two vectors in $\R^2$, we set $u\wedge v=x_1y_2-x_2y_1$. The following lemma is elementary, but will be useful for us in the sequel.

\begin{lemma}\label{SlitTorLm1}
Let $T$ be a flat torus, and $s$ be a geodesic segment joining two distinct points $x_1,x_2$ in $T$. Let $c$ be a simple closed
geodesic passing through $x_1$, not parallel to $s$. Then $x_1$ is the unique intersection point of $c$ and $s$ if and only if  $|V(s)\wedge V(c)| <\Aa(T)$.
\end{lemma}

\dem Using $\SL(2,\R)$, we can assume that $V(c)$ is horizontal and $V(s)$ is vertical. Cutting $T$ along $c$, we then get a cylinder $C$. Let
$h$ be the height of $C$. The fact that $x_1$ is the unique intersection point of $s$ and $c$ is equivalent to the fact that $|V(s)|<h$, which
is equivalent to

$$ |V(s)\wedge V(c)|=|V(s)||V(c)| < h|V(c)|=\Aa(C)=\Aa(T).$$

\carre

\subsection{The space of splittings}

Let $\Sig$ be a surface in $\H^\mathrm{hyp}(4)$. We denote by $P$ the unique singularity of $\Sig$. Let  $\delta_i^\pm, i=1,\dots,4,$ be four pairs of saddle connections in $\Sig$ as in Theorem \ref{ThA}. Cutting $\Sig$ along $(\delta_1^+,\delta_1^-)$ and $(\delta_3^+,\delta_3^-)$, we get three following components:

\begin{itemize}

\item[$\bullet$] $C_1$ is a cylinder bounded by $\delta_1^+$ and $\delta_1^-$. Gluing $\delta_1^+$ and $\delta_1^-$ together so that the two
points corresponding to $P$ are identified, we then get a surface in $\H(0)$ with a marked saddle connections.

\item[$\bullet$] $C_2$ is an annulus equipped with a flat metric structure with piecewise geodesic boundary, each boundary component of $C_2$ consists of two geodesic segments (corresponding to $\delta_1^+\cup\delta_3^+$, and $\delta_1^-\cup\delta_3^-$). Gluing $\delta_1^+$ and  $\delta_3^+$ to $\delta_1^-$ and $\delta_3^-$ respectively, we then get  an element of $\H(0,0)$, together with two saddle connections whose union is a simple closed curve.

\item[$\bullet$] $C_3$ is a one holed flat torus,  the boundary of $C_3$ is connected and consists of two geodesic segments
corresponding to $\delta_3^+$ and $\delta^-_3$. Gluing these two segments together, we then get an element in $\H(0,0)$ together with a marked
saddle connection.

\end{itemize}

\rem We get a similar decomposition of $\Sig$ by cutting along the pairs $(\delta_2^+,\delta_2^-)$ and $(\delta_4^+,\delta_4^-)$.\\

Let $\Sp$ denote the set of $(T_1,T_2,T_3, v_1,v_2)$, where $T_1 \in \H(0), T_2,T_3\in \H(0,0),$ and $v_i\in \R^2, i=1,2$ satisfying

\begin{itemize}

\item[a)] There are a saddle connection in $T_1$ and a saddle connection in $T_2$ both have associated vector equal to $v_1$.

\item[b)] There are a saddle connection in $T_2$, and a saddle connection in $T_3$ both have associated vector equal to $v_2$.

\item[c)] $v=v_1+v_2$ is a primitive vector of the lattice $\Lambda_2=\Lambda(T_2)$ associated to $T_2$, and there exists another primitive vector $w$ such that $\Lambda_2=\Z v\oplus\Z w$ and $0<|v_i\wedge w|<\Aa(T_2), i=1,2$.

\end{itemize}

\noindent We denote by $\Sp_1$ the subset of $\Sp$ consisting of elements $(T_1,T_2,T_3,v_1,v_2)$ such that $\Aa(T_1)+\Aa(T_2)+\Aa(T_3)=1$.\\

\rem

\begin{itemize}

\item[$\bullet$] We have a natural action of $\SL(2,\R)$ on $\Sp$.

\item[$\bullet$] It follows from the condition {\rm c)} that the flat torus $T_2$ is obtained from the gluing of two parallelograms
$\P_1,\P_2$ with $\P_i$ constructed from $w$ and $v_i$. 

%

\end{itemize}

\noindent Given an element $(T_1,T_2,T_3,v_1,v_2)$ in $\Sp$, we construct a surface in $\H^\mathrm{hyp}(4)$ as follows:

\begin{itemize}

\item[.] Cutting $T_1$ along the saddle connection corresponding to $v_1$, we get a cylinder $C_1$.

\item[.] Let $s_1,s_2$ be the saddle connections in $T_2$ corresponding to $v_1$ and $v_2$ respectively.  Since $v_1+v_2$ is a primitive vector in $\Lambda(T_2)$, we see that $s_1\cup s_2$ is a simple closed curve in $T_2$. Cutting $T_2$ along $s_1$ and $s_2$, we then get a cylinder with piecewise geodesic boundary, which will be denoted by $C_2$. 

\item[.] Slitting open $T_3$ along saddle connection corresponding to $v_2$, we get a one holed torus which will be denoted by $C_3$.

\item[.] We can now glue $C_1,C_2,C_3$ together following the model shown in Figure 1 so that all the marked points are identified, we then get a surface in $\H^\mathrm{hyp}(4)$.
\end{itemize}

\noindent This construction provides us with a map $\Psi: \Sp \lra \H^\mathrm{hyp}(4)$. A direct consequence of Theorem \ref{ThA} is the following

\begin{proposition}\label{ProjPr}
The map $\Psi$ is surjective, locally homeomorphic, and $\SL(2,\R)$-equivariant.\\
\end{proposition}

\subsection{Special splitting}\label{SpSplSect}

Let $X=\SpEmt$ be an element of $\Sp$, we say that $X$ is a {\em special splitting} if $v_1$ and $v_2$ are parallel (collinear). We denote by
$\SSp$ the set of special splittings in $\Sp$, and by $\SSp_1$ the intersection $\SSp\cap\Sp_1$.\\
Consider a point $X=\SpEmt$ in $\SSp$, we denote by $\Lambda_i, i=1,2,3,$ the lattices associated to $T_i$. Let $C_1$ (resp. $C_2$) denote the cylinder obtained by cutting $T_1$ (resp. $T_2$) along the saddle connection corresponding to $v_1$ (resp. along the union of the saddle connections corresponding to $v_1$ and $v_2$). Let $m_i, i=1,2,$ denote the modulus of $C_i$, we will call $m_1$ (resp. $m_2$) the modulus of the pair $(T_1,v_1)$ (resp. of the pair $(T_2,v_1+v_2)$). By construction, $C_1$ and $C_2$ are isometric to two cylinders in the direction $v_1$ on the surface $\Sig=\Psi(X)$. Set

$$\alpha=\frac{|v_2|}{|v_1|}, \text{ and } \bar{m}=\frac{m_1}{m_2}.$$

\noindent Observe that we have the following relation between $\bar{m}$ and $\alpha$

$$\bar{m}=\frac{\Aa(T_1)}{\Aa(T_2)}(1+\alpha)^2. $$

\noindent Since $\alpha$ and the areas of $T_i$ are $\SL(2,\R)$-invariant, so is  $\bar{m}$. We will call $\bar{m}$ the moduli ratio of
$\SpEmt$.\\
Using $SO(2,\R)$, we can assume that $C_1$ and $C_2$ are horizontal. We can also define the twists for $C_1$ and $C_2$ as follows: let $w_1=(w_1^x,w_1^y)$ (resp. $w_2=(w_2^x,w_2^y)$) be a primitive vector in $\Lambda_1$ (resp. $\Lambda_2$) such that $\Lambda_1=\Z v_1\oplus\Z w_1$ (resp. $\Lambda_2=\Z(v_1+v_2) \oplus\Z w_2$). We define the twists $t_1,t_2$ of $C_1$ and $C_2$ respectively to be $\DS{t_1=\frac{w_1^x}{|v_1|} \mod \Z, \text{ and } t_2=\frac{w_2^x}{|v_1|+|v_2|} \mod \Z}$. We also call $t_1$ (resp. $t_2$) the twist of the pair $(T_1,v_1)$ (resp. of the pair $(T_2,v_1+v_2)$).\\
Recall that a vector $w$ in $\R^2$ is generic with respect to a lattice $\Lambda= \Z u \oplus \Z v$ if $w$ is not parallel to any vector in $\Lambda$. To prove Theorem \ref{ThB}, we first prove the following theorem, which is slightly weaker. As we will see, Theorem \ref{ThB} can be obtained as a consequence of this theorem.

\begin{theorem}\label{ThC}
Let $X_0=(T^0_1,T^0_2,T^0_3,v^0_1,v^0_2)$  be an element in $\SSp_1$. Let $\Lambda^0_i, i=1,2,3,$ denote the lattice associated to $T^0_i$, and
$\bar{m}_0$ denote the moduli ratio of $X_0$. Suppose that

\begin{itemize}

\item[$\bullet$] $\bar{m}_0 \notin \Q$,

\item[$\bullet$] $v^0_2$ is generic with respect to $\Lambda^0_3$,
\end{itemize}

\noindent then $\Orb:=\SL(2,\R)\cdot \Psi(X_0)$ is dense in $\H^{\mathrm{hyp}}_1(4)$.

\end{theorem}


\subsection{Ratner's Theorem}

The first important ingredient of the proof of Theorem \ref{ThC} is a consequence of the famous theorem of Ratner on action of unipotent subgroups
on homogeneous spaces. Before stating this theorem, let us first recall some basic notions. Let $G$ be a Lie group, and $\mathfrak{g}$ be its
Lie algebra. An element $g$ of $G$ is {\em unipotent} if $\mathbf{Ad}_g -\Id$ is nilpotent in $\mathrm{End}(\mathfrak{g})$. Let $\lambda$ denote the right Haar measure of $G$, $G$ is called {\em unimodular} if the left Haar measure equals the right Haar measure, or equivalently if $|\det \mathbf{Ad}_g|=1$ for all $g$ in $G$. A discrete subgroup $\Gamma$ of $G$ is called a {\em lattice} if we have $\DS{\lambda(G/\Gamma)<\infty}$. If $G$ has a lattice then it is unimodular. It is well-known that $\SL(2,\R)$ is unimodular, but its subgroup consisting upper triangular matrices is not.

\begin{theorem}[Ratner]\label{RatnerTh}
Let $G$ be a finite dimensional Lie group, $\Gamma$ be a lattice in $G$, and $X=G/\Gamma$. Let $U$ be a connected subgroup of $G$ generated by
unipotent element. Then for any $x$ in $X$, the closure $\overline{U\cdot x}$ of the $U$-orbit of $x$ is a homogeneous space of finite volume,
that is there exists a closed unimodular subgroup $H \subset G$ containing $U$ such that

\begin{itemize}
\item[$\bullet$] $\overline{U\cdot x} = H\cdot x$,

\item[$\bullet$] $x\Gamma x^{-1}\cap H$ is a lattice in $H$.
\end{itemize}

\end{theorem}

Put $G=\R \times \R \times \SL(2,\R), \text{ and } \Gamma= \Z\times\Z\times\SL(2,\Z)$, then $\Gamma$ is a lattice in $G$. An element of $G/\Gamma$ is a triple $(\theta_1,\theta_2,\Lambda)$, where $\theta_i \in \R/\Z \simeq \S^1$, and $\Lambda \simeq \Z^2$ is a lattice in $\R^2$ such that $\mathrm{Vol}(\R^2/\Lambda)=1$. Let $m_1,m_2$ be two positive real numbers. We set

$$U=U_{m_1,m_2}=\{(m_1t, m_2t, \left(%
\begin{array}{cc}
  1 & t \\
  0 & 1 \\
\end{array}%
\right)), \; \; t\in \R \},$$

\noindent then $U$ is a unipotent subgroup of $G$. As a consequence of Theorem \ref{RatnerTh}, we have

\begin{corollary}\label{RatnerCor}
Suppose that $m_1/m_2 \notin \Q$. Let $\Lambda$ be a lattice in $\R^2$ which contains no horizontal vectors. Then for any
$(\theta_1,\theta_2)\in \R/\Z\times\R/\Z$, we have

$$\overline{U\cdot (\theta_1,\theta_2,\Lambda)}=G/\Gamma.$$

\end{corollary}

\dem By Ratner Theorem, we know that $\overline{U\cdot (\theta_1,\theta_2,\Lambda)}=H\cdot (\theta_1,\theta_2,\Lambda)$, where $H$ is connected,
unimodular subgroup of $G$. All we need to show is that $H=G$.\\
Let $x$ be any element of $G$ which is projected to $(\theta_1,\theta_2,\Lambda)$. Let $\mathfrak{h}$ and $\mathfrak{g}$ denote the
Lie algebras of $H$ and $G$ respectively. Set

$$\mathbf{a}=\left(%
\begin{array}{cc}
  1 & 0 \\
  0 & -1 \\
\end{array}%
\right), \mathbf{u}_+=\left(%
\begin{array}{cc}
  0 & 1 \\
  0 & 0 \\
\end{array}%
\right), \mathbf{u}_-=\left(%
\begin{array}{cc}
  0 & 0 \\
  1 & 0 \\
\end{array}%
\right).$$

\noindent We have $\mathfrak{sl}(2,\R)=\R\mathbf{u}_+\oplus\R\mathbf{u}_-\oplus\R\mathbf{a}, \text{ and  } \mathfrak{g}=\R\oplus\R\oplus\mathfrak{sl}(2,\R)$. Observe that the Lie bracket of $\mathfrak{g}$ is trivial on the $\R$ components, and we have

$$[\mathbf{u}_+,\mathbf{u}_-]=\mathbf{a}, [\mathbf{a},\mathbf{u}_+] =2\mathbf{u}_+, [\mathbf{a},\mathbf{u}_-]=-2\mathbf{u}_-.$$

\noindent Since $U\subset H$, the Lie algebra $\mathfrak{h}$ contains $\mathbf{v}_0=(m_1,m_2,\mathbf{u}_+)$. Remark that $x\Gamma
x^{-1}=\Z\times\Z\times M\cdot\SL(2,\Z)\cdot M^{-1}$, where $M$ is any matrix in $\SL(2,\R)$ sending the standard basis of $\R^2$ to a basis of
the lattice $\Lambda$. We denote by $A$ and $N$ the following subgroups of $\SL(2,\R)$

$$ A=\{\left(%
\begin{array}{cc}
  e^t & 0 \\
  0 & e^{-t} \\
\end{array}%
\right), \; t\in \R\}, \; N=\{\left(%
\begin{array}{cc}
  1 & t \\
  0 & 1 \\
\end{array}%
\right),\; t\in \R\}.$$

\noindent Let $\mathbf{pr}_2: \mathfrak{g} \lra \mathfrak{sl}(2,\R)$ denote the natural projection. The image of $\mathfrak{h}$ under
$\mathbf{pr}_2$ is a subalgebra of $\mathfrak{sl}(2,\R)$ which contains $\R\mathbf{u}_+$.\\

\noindent\underline{\bf Case 1:} $\mathbf{pr}_2(\mathfrak{h})=\R\mathbf{u}_+$. We have three possibilities:

\begin{itemize}
\item[$\bullet$] $\mathfrak{h}=\R\mathbf{v}_0 \Longrightarrow H=U$, but by assumption, we have $U\cap \Z\times\Z\times M\cdot \SL(2,\Z) \cdot M^{-1}=\{(0,0,\Id)\}$ is not a lattice in $U$.

\item[$\bullet$] $\mathfrak{h}=\R \mathbf{v}_0\oplus \R\mathbf{u}_+=\R\mathbf{u}_+\oplus \R_{m_1,m_2}$, where $\R_{m_1,m_2}=\R\cdot(m_1,m_2) \subset \R^2$. It follows that $H=\R_{m_1,m_2}\times N$. But again, we have $H\cap \Z\times\Z\times M \cdot\SL(2,\Z)\cdot
M^{-1}=\{(0,0,\Id)\}$.

\item[$\bullet$] $\mathfrak{h}=\R^2\oplus \R\mathbf{u}_+ \Longrightarrow H=\R^2\times N$. But we have $N\cap M\cdot\SL(2,\Z)\cdot
M^{-1}=\{\Id\}$, therefore, $H\cap \Z\times\Z \times M\cdot\SL(2,\Z)\cdot M^{-1}$ is not a lattice.\\
\end{itemize}

\noindent \underline{\bf Case 2:} $\mathbf{pr}_2(\mathfrak{h})=\R\mathbf{u}_+\oplus\R\mathbf{a}$. Let $v$ be any vector in $\mathfrak{h}$ such
that $\mathbf{pr}_2(\mathbf{v})=\mathbf{a}$, then we have $[\mathbf{v},\mathbf{v}_0]=2\mathbf{u}_+$. Therefore, we see that $\mathfrak{h}$ contains the following vectors

\begin{itemize}
\item[.] $\mathbf{u}_+$,

\item[.] $\mathbf{v}_1=\mathbf{a}+\mathbf{w}$, with $\mathbf{w}=(k_1,k_2)\in \R^2$,

\item[.] $\mathbf{w}_0=\mathbf{v_0}-\mathbf{u}_+=(m_1,m_2) \in \R^2$.

\end{itemize}

\noindent Here we have two possibilities:

\begin{itemize}

\item[$\bullet$] $\mathfrak{h}=\R\mathbf{w}_0\oplus\R\mathbf{u}_+\oplus\R\mathbf{v}_1 \Longrightarrow H =\R\times A'N$, where $A'=\{(k_1t,k_2t,\left(%
\begin{array}{cc}
  e^t & 0 \\
  0 & e^{-t} \\
\end{array}%
\right)), \; t\in \R \} \subset G$. It follows that $H\simeq \R\times AN.$ But since $AN$ is not unimodular, neither is $H$.\\

\item[$\bullet$]  $\mathfrak{h}=\R^2\oplus\R\mathbf{u}_+\oplus\R\mathbf{a} \Longrightarrow H =\R^2\times AN$, but again $H$ is not unimodular.

\end{itemize}

\noindent \underline{\bf Case 3:} $\mathbf{pr}_2(\mathfrak{h})=\mathfrak{sl}(2,\R)$. Let $\mathbf{v}$ be a vector in $\mathfrak{h}$ such that
$\mathbf{pr}_2(\mathbf{v})=\mathbf{u}_-$, we then have

$$[\mathbf{v}_0,\mathbf{v}]= \mathbf{a}, \; [[\mathbf{v}_0,\mathbf{v}],\mathbf{v_0}]=2\mathbf{u}_+ \text{ and } [[\mathbf{v}_0,\mathbf{v}],\mathbf{v}]=-2\mathbf{u}_-.$$

\noindent It follows  that $\mathfrak{sl}(2,\R)\subset \mathfrak{h}$. We then have two possibilities:

\begin{itemize}
\item[$\bullet$] $H=\R_{m_1,m_2}\times \SL(2,\R)$, in this case $H\cap \Z\times\Z\times M\cdot\SL(2,\Z)\cdot M^{-1}=(0,0,M\cdot\SL(2,\Z)\cdot
M^{-1})$ is not a lattice in $H$.

\item[$\bullet$] $H=\R^2\times\SL(2,\R)$, this is the only admissible possibility.
\end{itemize}

We can then conclude that $\overline{U\cdot(\theta_1,\theta_2,\Lambda)}=G/\Gamma$. \carre

\rem  Similar results for $\R^k\times\SL(2,\R)^n/\Z^k\times\SL(2,\Z)^n$ with small $k$ and $n$ can be found in \cite{HubLanMo3}.\\

For any $(A_1,A_2,A_3,\alpha)$ in $\R_{>0}^4$, let $\SSp(A_1,A_2,A_3,\alpha)$ denote the subset of $\SSp$ consisting of elements $(T_1,T_2,T_3,v_1,v_2)$ such that $\Aa(T_i)=A_i, i=1,2,3,$ and $\DS{\frac{|v_2|}{|v_1|}=\alpha}$. Using Corollary \ref{RatnerCor}, we have the following lemma

\begin{lemma}\label{PrfCLm1}

Let $X=(T_1,T_2,T_3,v_1,v_2)$ be an element in $\SSp(A_1,A_2,A_3,\alpha)$. If $\DS{\frac{A_1}{A_2}(\alpha+1)^2}$ is irrational, and $v_2$ is generic with respect to the lattice $\Lambda_3=\Lambda(T_3)$ then 

$$\Psi(\SSp(A_1,A_2,A_3,\alpha)) \subset \overline{\SL(2,\R)\cdot\Psi(X)}.$$

\end{lemma}

\dem   Let $\SSp(A_1,A_2,A_3,\alpha)_{\mathrm{hor}}$ denote the subset of $\SSp(A_1,A_2,A_3,\alpha)$ consisting of elements with $v_1=(1,0)$. We have $\SSp(A_1,A_2,A_3,\alpha)=\SL(2,\R)\cdot\SSp(A_1,A_2,A_3,\alpha)_\mathrm{hor}$. We have a natural mapping $\varphi :\SSp(A_1,A_2,A_3,\alpha)_\mathrm{hor} \lra G/\Gamma$ which sends $(T_1,T_2,T_3,v_1,v_2)$ to an element $(t_1,t_2,\Lambda_3)$, where $t_1$ and $t_2$ are the twists of $(T_1,v_1)$ and $(T_2,v_1+v_2)$ respectively, and $\Lambda_3$ is the lattice associated to $T_3$ normalized to have covolume one. Remark that $\varphi$ is a homeomorphism onto its image.\\
Let $m_1$ and $m_2$ denote the moduli of $(T_1,v_1)$ and $(T_2,v_1+v_2)$ respectively. Recall that we have $\DS{\frac{m_1}{m_2}=\frac{A_1}{A_2}(\alpha+1)^2}$. We define the action of $U=\{\left(%
\begin{array}{cc}
  1 & t \\
  0 & 1 \\
\end{array}%
\right), \; t\in\R \}$ on $G/\Gamma$ using the identification $U\simeq U_{m_1,m_2}$. It follows that $\varphi$ is $U$-equivariant. \\
Without loss of generality, we can assume that $X\in \SSp(A_1,A_2,A_3,\alpha)_\mathrm{hor}$. Consider $x=\varphi(X)\in G/\Gamma$.
The hypothesis on $X$ implies that $x$ satisfies the conditions of Corollary \ref{RatnerCor}, therefore $\overline{U\cdot x}=G/\Gamma$.
Since $\varphi$ is $U$-equivariant and a local homeomorphism, we deduce that $\overline{U\cdot X}=\varphi^{-1}(\overline{U\cdot x})=\SSp(A_1,A_2,A_3,\alpha)_\mathrm{hor}$, and the lemma follows. \carre

\begin{corollary}\label{PrfCCor1}
Let $X_0=(T_1^0,T_2^0,T_3^0,v_1^0,v_2^0)$ be as in Theorem \ref{ThC}. Then we have 
$$\Psi(\SSp(A_1^0,A_2^0,A_3^0,\alpha_0)) \subset\overline{\SL(2,\R)\cdot\Psi(X_0)},$$

\noindent where $A_i^0=\Aa(T_i^0), i=1,2,3$, and $\alpha_0=|v_2^0|/|v_1^0|$.
\end{corollary}

\section{Surfaces admitting special splitting are contained in the orbit closure}

Our aim in this section is to prove the following

\begin{proposition}\label{PropC}

Let $X_0=(T_1^0,T_2^0,T_3^0, v_1^0,v_2^0)$ be as in Theorem \ref{ThC}. We have $\Psi(\SSp_1)\subset \overline{\SL(2,\R)\cdot \Psi(X_0)}$.

\end{proposition}

\subsection{Dual splitting}

Given $X=(T_1,T_2,T_3,v_1,v_2)$ in $\SSp$, we will denote both saddle connections in $T_1$ and $T_2$ corresponding to $v_1$ by $\delta_1$, similarly, we denote by $\delta_2$ the two saddle connections in $T_2$ and $T_3$ corresponding to $v_2$. Recall that the saddle connections $\delta_i, i=1,2,$ give rises to a pair of homologous saddle connections in the surface $\Sig=\Psi(X)$, which will be denoted by $\delta_i^\pm$.

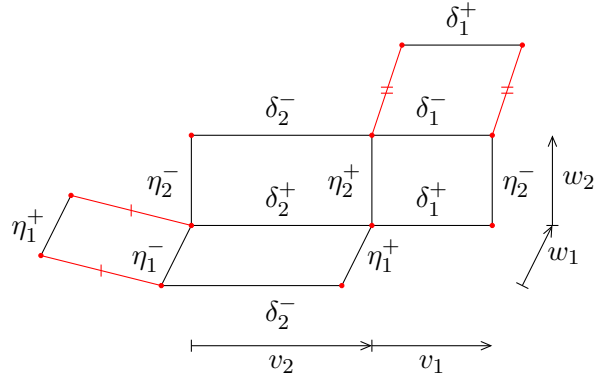
\begin{figure}[!ht]\label{Fig4}
\begin{center}

\begin{tikzpicture}[scale=0.4]

\draw (-1,2) -- (0,4)  (-4,5) -- (-5,3)  (-1,2) -- (5,2) -- (6,4) -- (10,4) -- (10,7)  (11,10) -- (7,10)  (6,7) -- (0,7) -- (0,4) -- (6,4) -- (6,7) -- (10,7);

\draw[red] (0,4) -- (-4,5) (-5,3) -- (-1,2) (10,7) -- (11,10) (7,10) -- (6,7);

\foreach \x in {(-2,4.5),(-3,2.5)} \draw[red] \x +(0,-0.2) -- +(0,0.2);
\foreach \x in {(10.5, 8.5),(6.5, 8.5)} \draw[red] \x +(-0.2,0.1) -- +(0.2,0.1) +(-0.2,-0.1) -- +(0.2,-0.1);

\draw[>=angle 45, ->] (0,0) -- (6,0); \draw[>=angle 45, ->] (6,0) -- (10,0); \draw (0,0) +(0,-0.2) -- +(0,0.2) (6,0) +(0,-0.2) -- +(0,0.2);

\foreach \x in {(-1,2), (5,2), (6,4), (10,4), (10,7), (11,10), (7,10), (6,7), (0,7), (0,4), (-4,5), (-5,3)} \filldraw[red] \x circle (2pt);

\draw (9,10) node[above] {$\delta_1^+$} (8,7) node[above] {$\delta_1^-$} (8,4) node[above] {$\delta_1^+$} (3,7) node[above] {$\delta_2^-$} (3,4) node[above] {$\delta_2^+$} (3,2) node[below] {$\delta_2^-$};

\draw (-4.5,4) node[left] {$\eta_1^+$} (-0.5,3) node[left] {$\eta_1^-$} (5.5,3) node[right] {$\eta_1^+$} (0,5.5) node[left] {$\eta_2^-$} (6,5.5) node[left] {$\eta_2^+$} (10,5.5) node[right] {$\eta_2^-$};

\draw (3,0) node[below] {$v_2$} (8,0) node[below] {$v_1$};

\draw[>=angle 45, ->] (11,2) -- (12,4); \draw[>=angle 45, ->] (12,4) -- (12,7); \draw (11,2) +(-0.2,0.1) -- +(0.2,-0.1) (12,4) +(-0.2,0) -- +(0.2,0);

\draw (11.5,3) node[right] {$w_1$} (12,5.5) node[right] {$w_2$};

\end{tikzpicture}

\caption{Dual splittings}

\end{center}

\end{figure}

\noindent Let $\eta_1$ be a simple closed geodesic in $T_3$ which meets $\delta_2$ once, and let $\eta_1^{\pm}$ denote the pair of saddle connections parallel to $\eta_1$. Similarly, let $\eta_2$ be a simple closed geodesics in $T_2$ which meets $\delta_1\cup\delta_2$ once, and let $\eta_2^{\pm}$ denote the pair of saddle connections parallel to $\eta_2$. Remark that $\eta_1^{\pm}$ (resp. $\eta_2^{\pm}$) are homologous saddle connections in $\Sig$. We choose the orientation of $\eta_1$ and $\eta_2$ so that $\eta^+_1*\eta^+_2$ is freely homotopic to a simple closed curve in $\Sig$. Put $w_1=V(\eta_1), w_2=V(\eta_2)$.\\
Cutting $\Sig$ along $\eta_1^{\pm}$ and $\eta_2^{\pm}$, we see that the surface $\Sig$ is obtained from another element $X^\vee=(T^\vee_1,T^\vee_2,T^\vee_3,w_1,w_2)$ in $\Sp$.  We will call $X^\vee$ a {\em dual splitting} of $X$. Note that $X^\vee$ does not belong to $\SSp$ in general, and there are infinitely many splittings dual to a given splitting. We also have

\begin{equation}\label{DualArea}
\Aa(T_3^\vee)=\Aa(T_1)+\frac{\Aa(T_2)}{1+|v_2|/|v_1|}
\end{equation}

\noindent Throughout this section, we set $A_i^0=\Aa(T_i^0), i=1,2,3$, and $\DS{\alpha_0=\frac{|v_2^0|}{|v_1^0|}}$.

\subsection{ Changing splitting}
The first step to prove Proposition \ref{PropC} is the following

\begin{lemma}\label{PrfCLm2}
If $(A_1,A_2,A_3,\alpha)\in \R_{>0}^4$ satisfies

\begin{itemize}
\item[.] $A_1+A_2+A_3=1$,

\item[.] $\DS{A_1+\frac{A_2}{1+\alpha}=A_1^0+\frac{A_2^0}{1+\alpha_0}}$
\end{itemize}

\noindent then $\Psi(\SSp(A_1,A_2,A_3,\alpha))\subset \overline{\SL(2,\R)\cdot\Psi(X_0)}$.
\end{lemma}

By Corollary \ref{PrfCCor1}, we know that $\overline{\SL(2,\R)\cdot \Psi(X_0)}$ contains $\Psi(\SSp(A_1^0,A_2^0,A_3^0,\alpha_0))$. Let
$X=(T_1,T_2,T_3,v_1,v_2)$ be an element in $\SSp(A_1^0,A_2^0,A_3^0,\alpha_0)$, and $\Sig$ be the surface in $\H^\mathrm{hyp}(4)$ constructed
from $X$. Let $\delta_i,\delta_i^\pm, \eta_i, \eta_i^\pm, i=1,2$, and $X^\vee=(T_1^\vee,T_2^\vee,T_3^\vee,w_1,w_2)$ be as in the previous subsection, where $X^\vee$ is a dual splitting of $X$.\\
Let $\sigma^\pm_1$ (resp. $\sigma_2^\pm$) be a pair of homologous saddle connections in $T^\vee_3$ (resp. $T^\vee_2$) which bound a simple cylinder containing $\eta_2$ (see Figure 5). Viewed as saddle connections of $\Sig$,  the pairs $\sigma_1^\pm$ and $\sigma_2^\pm$  determine a splitting of $\Sig$. If $\sigma_1^\pm$ and $\sigma_2^\pm$ are parallel, then we have another special splitting of $\Sig$. To prove the lemma, we will show that for any $(A_1,A_2,A_3,\alpha)$ in $\R_{>0}^4$, there exists an element $X$ in $\SSp(A_1^0,A_2^0,A_3^0,\alpha_0)$ for which one can find $\sigma_1^\pm, \sigma_2^\pm$ determining a special splitting with parameters $(A_1,A_2,A_3,\alpha)$. We can then use Lemma \ref{PrfCLm1} to conclude, first, for $(A_1,A_2,A_3,\alpha)$ satisfying the condition of Lemma \ref{PrfCLm1}, and then for all $(A_1,A_2,A_3,\alpha)$ by continuity.\\

\dem (of Lemma \ref{PrfCLm2}) Without loss of generality, we can assume that $v_1=(1,0)$ and $v_2=(\alpha_0,0)$. Let $C_1$ (resp. $C_2$)
denote the cylinder obtained by slitting $T_1$ (resp. $T_2$) along saddle connection $\delta_1$ (resp. along the saddle connections $\delta_1$ and $\delta_2$). Let $h_i$ and $t_i$ denote the height and the twist of $C_i, i=1,2$. Note that $h_1=A_1^0$, and $\DS{h_2=\frac{A_2^0}{\alpha_0+1}}$. 
We fix $t_2=0$, consequently, we can choose  $\eta_2^\pm$ to be vertical, and therefore $w_2=V(\eta_2^\pm)=(0,h_2)$.\\
Set $\Lambda_i=\Lambda(T_i)$ and $\Lambda^\vee_i=\Lambda(T^\vee_i), i=1,2,3$. Recall that $v_1$ is a primitive vector of $\Lambda_1$, let $u_1=(x,h_1)$ be another primitive vector in $\Lambda_1$ such that $\Lambda_1=\Z u_1\oplus \Z v_1$. Observe that the parameter $x$ can be chosen arbitrarily. Similarly, $w_1=(y,z)$ is a primitive vector in $\Lambda_1^\vee$, let $\hat{u}_1$ be another primitive vector such that $\Lambda^\vee_1=\Z w_1 \oplus \Z \hat{u}_1$. Note that $\Lambda_2^\vee=\Z v_2 \oplus \Z(w_1+w_2)$ and $\Lambda_3=\Z w_1 \oplus \Z (v_2+\hat{u}_1)$.
The parameters $(x,y,z,\hat{u}_1)\in \R^3\times \R^2$ uniquely determine the element $X$ in $\SSp(A_1^0,A_2^0,A_3^0,\alpha_0)$. By construction, the parameters $(x,y,z, \hat{u}_1)$ must satisfy the following conditions

\begin{equation}\label{PrfCIneq1}
|v_2\wedge w_2| < \Aa(T^\vee_2)=|v_2\wedge(w_1+w_2)| < 1-\Aa(T^\vee_3)=\frac{\alpha_0A^0_2}{\alpha_0+1}+A^0_3
\end{equation}

\begin{equation}\label{PrfCIneq2}
|w_1\wedge (v_2+\hat{u}_1)|=A_3^0 
\end{equation}

\noindent Simple computations show that (\ref{PrfCIneq1}) is equivalent to

\begin{equation}\label{CondZ}
0< z < \frac{A_3^0}{\alpha_0}
\end{equation}

\noindent Remark that the conditions (\ref{PrfCIneq1}) and (\ref{PrfCIneq2}) are sufficient, that is, if the parameters $(x,y,z,\hat{u}_1)$ satisfy these two conditions, then they determine an element in  $\SSp(A_1^0,A_2^0,A_3^0,\alpha_0)$.

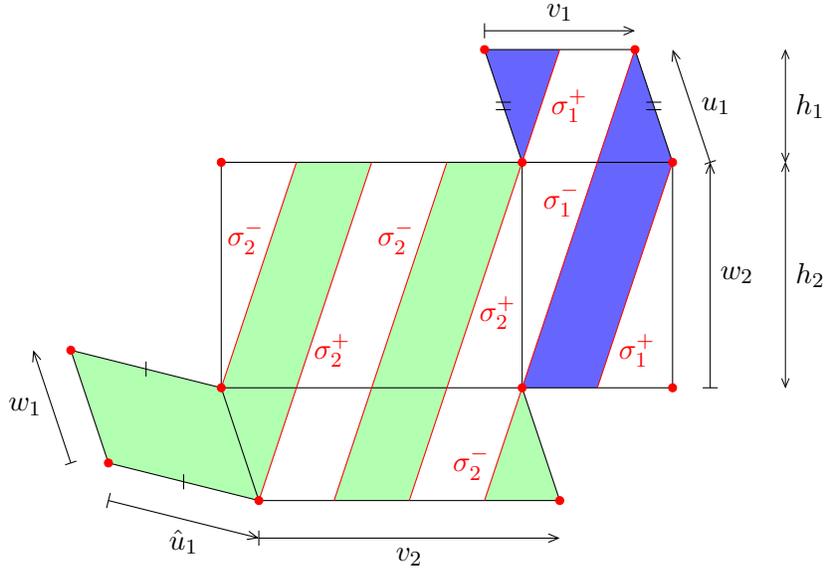
\begin{figure}[!ht]\label{Fig5}

\begin{center}
\begin{tikzpicture}[scale=0.5]

\fill[blue!60] (8,9) -- (9,12) -- (7,12) -- cycle; 
\fill[blue!60] (8,3) -- (11,12) -- (12,9) -- (10,3) -- cycle;

\fill[green!30] (8,3) -- (9,0) -- (7,0)-- cycle;
\fill[green!30] (5,0) -- (8,9) -- (6,9) -- (3,0) -- cycle;
\fill[green!30] (1,0) -- (4,9) -- (2,9) -- (0,3) -- (-4,4) -- (-3,1) -- cycle;

\draw (1,0) -- (0,3) -- (-4,4) -- (-3,1) -- (1,0) -- (9,0) -- (8,3) -- (12,3) -- (12,9) -- (11,12) -- (7,12) -- (8,9) -- (0,9) -- (0,3) -- (8,3) -- (8,9) -- (12,9);

\draw (-2,3.5) +(0,-0.2) -- +(0,0.2) (-1,0.5) +(0,-0.2) -- +(0,0.2);
\draw (7.5,10.5) +(-0.2,-0.1) -- +(0.2,-0.1) +(-0.2,0.1) -- +(0.2,0.1) (11.5,10.5) +(-0.2,-0.1) -- +(0.2,-0.1) +(-0.2,0.1) -- +(0.2,0.1);

\draw[red] (0,3) -- (2,9) (1,0) -- (4,9) (3,0) -- (6,9) (5,0) -- (9,12)  (7,0) -- (11,12) (10,3) -- (12,9);

\foreach \x in {(1,0),(0,3),(-4,4),(-3,1),(9,0),(8,3),(12,3), (12,9),(11,12),(7,12),(8,9),(0,9)} \filldraw[red] \x circle (3pt);

\draw[red] (8.5,10.5) node[right] {$\sigma_1^+$} (9.8,8) node[left] {$\sigma_1^-$} (10.3,4) node[right] {$\sigma_1^+$}; 

\draw[red] (6.6,5) node[right] {$\sigma_2^+$} (7.4,1) node[left] {$\sigma_2^-$} (5.4,7) node[left] {$\sigma_2^-$} (2.2,4) node[right] {$\sigma_2^+$} (1.4,7) node[left] {$\sigma_2^-$};

\draw[>=angle 45, ->] (7,12.5) -- (11,12.5); \draw (7,12.5) +(0,-0.2) -- +(0,0.2); \draw (9,12.5) node[above] {$v_1$};

\draw[>= angle 45, ->] (13,9) -- (12,12); \draw (13,9) +(-0.15, -0.05) -- +(0.15,0.05); \draw (12.5,10.5) node[right] {$u_1$};

\draw[>= angle 45, ->] (-4,1) -- (-5,4); \draw (-4,1) +(-0.15, -0.05) -- +(0.15,0.05); \draw (-4.5,2.5) node[left] {$w_1$};

\draw[>= angle 45, <-] (1,-1) -- (-3,0); \draw (-3,0) +(-0.05,-0.2) -- +(0.05,0.2); \draw (-1,-0.5) node[below] {$\hat{u}_1$};

\draw[>=angle 45, ->]  (1,-1) -- (9,-1); \draw (1,-1) +(0,-0.2) -- +(0,0.2); \draw (5,-1) node[below] {$v_2$};

\draw[>=angle 45, ->] (13,3) -- (13,9); \draw (13,3) +(-0.2,0) -- +(0.2,0); \draw (13,6) node[right] {$w_2$};

\draw[>= angle 45, <->, very thin] (15,9) -- (15,12); \draw[>= angle 45, <->, very thin] (15,9) -- (15,3);

\draw (15,10.5) node[right] {$h_1$} (15,6) node[right] {$h_2$};

\end{tikzpicture}
\caption{Finding new special splittings}
\end{center}

\end{figure}

\noindent \underline{\bf Claim 1:} \textit{ For any $(A_1,A_2,A_3,\alpha) \in \R_{>0}^4$ satisfying the conditions of the lemma, there exist $(x,y,z)\in \R^3$ with $z$ satisfying (\ref{CondZ}) such that we can find a primitive vector $v'_1$ in $\Lambda^\vee_3$, and a primitive vector $v'_2$ in $\Lambda^\vee_2$ such that}

\begin{itemize}

\item[i)] $v'_2=\alpha v'_1$,

\item[ii)]  $\DS{|v'_1\wedge w_2|=\frac{A_2}{\alpha+1}}$,

\item[iii)] $|v'_2\wedge w_2| < \Aa(T^\vee_2)$.

\end{itemize}

\noindent \underline{\bf Proof of Claim 1:} Recall that, by  assumption, we have $\DS{A_1+\frac{A_2}{\alpha+1}=A_1^0+\frac{A_2^0}{\alpha_0+1}}$. Since $A_1+A_2<A_1^0+A_2^0+A_3^0=1$, it follows

\begin{equation}\label{PrfCIneq3}
\frac{\alpha A_2}{\alpha+1}<\frac{\alpha_0A_2^0}{\alpha_0+1}+A_3^0
\end{equation}

\noindent From (\ref{PrfCIneq3}), we deduce that there exist $p,q \in \N, p>0, q>0,$ such that

$$
\max\{\frac{\alpha_0A_2^0/(\alpha_0+1)}{\alpha\alpha_0(h_1+h_2)}, \frac{\alpha A_2/(\alpha+1)}{\alpha\alpha_0(h_1+h_2)}\} < \frac{p}{q}  <
\frac{\alpha_0A_2^0/(\alpha_0+1)+A_3^0}{\alpha\alpha_0(h_1+h_2)}$$

\noindent Set $\DS{x=\frac{1}{p}(\frac{A_2}{h_2(\alpha+1)}-1), y=\frac{1}{q}(\frac{\alpha A_2}{h_2(\alpha+1)}-\alpha_0), z=\alpha\frac{p}{q}(h_1+h_2)-h_2}$. By the choice of $p,q,$ it is straight forward to verify that $z$ satisfies (\ref{CondZ}). We have

\begin{eqnarray*}
u_1 & = & (x,h_1) = (\frac{A_2}{ph_2(\alpha+1)}-\frac{1}{p}, h_1)\\
w_1 & = & (y,z) = (\frac{\alpha A_2}{qh_2(\alpha+1)}-\frac{\alpha_0}{q}, \alpha\frac{p}{q}(h_1+h_2)-h_2)
\end{eqnarray*}

\noindent Set
\begin{eqnarray*}
v'_1 & = & v_1+p(u_1+w_2)=(1,0)+(\frac{A_2}{h_2(\alpha+1)}-1, p(h_1+h_2))=(\frac{A_2}{h_2(\alpha+1)}, p(h_1+h_2))\\
v'_2 & = & v_2+q(w_1+w_2)=(\alpha_0,0)+(\frac{\alpha A_2}{h_2(\alpha+1)}-\alpha_0, \alpha p(h_1+h_2))=\alpha(\frac{A_2}{h_2(\alpha+1)},p(h_1+h_2))
\end{eqnarray*}

\noindent Since $\Lambda^\vee_3$ is generated by $v_1$ and $u_1+w_2$, we see that $v'_1$ is a primitive vector in $\Lambda^\vee_3$, similarly, $v'_2$ is a primitive vector in $\Lambda_2^\vee$. Clearly, we have $v'_2=\alpha v'_1$, hence {\rm i)} is satisfied. We have 

$$|v'_1\wedge w_2|= \left|\begin{array}{cc}
  A_2/(h_2(\alpha+1)) & 0 \\
  p(h_1+h_2) & h_2 \\
\end{array}\right|= \frac{A_2}{\alpha+1}$$

\noindent therefore {\rm ii)} is satisfied. Next, we have 
$$\DS{\Aa(T^\vee_2)=|v_2\wedge (w_1+w_2)|=\alpha\alpha_0\frac{p}{q}(h_1+h_2)},$$ 
\noindent and  
$$\DS{|v'_2\wedge w_2|=\alpha|v'_1\wedge w_2|=\frac{\alpha A_2}{\alpha+1}}.$$

\noindent  By the choice of $p,q$, we have  $|v'_2\wedge w_2|<\Aa(T^\vee_2)$,  hence  {\rm iii)} is satisfied. \carre

\noindent \underline{\bf Claim 2:} \textit{Given  $(x,y,z)$ as in Claim 1, there exist $\hat{u}_1$ satisfying (\ref{PrfCIneq2}) such that $v'_2$ is generic with respect to the lattice $\Z w_1\oplus \Z(v'_2+\hat{u}_1)$.}\\

\noindent \underline{\bf Proof of Claim 2:}  Since $\DS{|v'_2\wedge w_1|=\alpha_0(h_1+h_2)\frac{p}{q}- \frac{A_2}{\alpha+1}>0}$, we deduce that $\{v'_2,w_1\}$ is a basis of $\R^2$. Therefore, we can write $\hat{u}_1=\lambda w_1+ \mu v'_2$. Observe that, once $w_1$ is fixed, the set of $\hat{u}_1$ satisfying (\ref{PrfCIneq2}) is parameterized by $\lambda\in\R$, with fixed $\mu$.\\
Observe that $v'_2$ is parallel to a vector in $\Z w_1+\Z(v'_2+\hat{u}_1)$ if and only if $\lambda \in \Q$. Indeed, suppose that $v'_2=\lambda'(mw_1+n(v'_2+\hat{u}_1))$, with $m,n\in \Z$, then we must have $n\neq 0$, otherwise $v'_2$ and $w_1$ are collinear, therefore $\DS{\hat{u}_1=-\frac{m}{n}w_1+\lambda''v'_2}$. It follows immediately that there exist $\hat{u}_1$ satisfying (\ref{PrfCIneq2}) such that $v'_2$ is generic with respect to $\Z w_1\oplus \Z(v'_2+\hat{u}_1)$. \carre

Let us now show that the lemma will follow from Claim 1 and Claim 2. Choose $(x,y,z)$ as in Claim 1, and choose $\hat{u}_1$ as in Claim 2, then the parameters $(x,y,z,\hat{u}_1)$ give us an element $X$ in $\SSp(A_1^0,A_2^0,A_3^0,\alpha_0)$. We have $\Sig=\Psi(X) \in \overline{\SL(2,\R)\cdot \Psi(X_0)}$.\\
Let $\sigma_1^\pm$ (resp. $\sigma_2^\pm$) be the pair of saddle connections in $T^\vee_3$ (resp. in $T^\vee_2$) corresponding to $v'_1$ (resp. $v'_2$). Since $\DS{|v'_1\wedge w_2|=\frac{A_2}{\alpha+1}<\Aa(T^\vee_3)=\frac{A^0_2}{\alpha_0+1}+A^0_1}$, from Lemma \ref{SlitTorLm1}, we deduce that $\sigma^\pm_1$ meet $\eta_2$ at only one point. Consequently, we see that $\sigma_1^\pm$ bound a simple cylinder containing $\eta_2$. Similarly, since $\Aa(T^\vee_2)=|v'_2\wedge w_1| + |v'_2\wedge w_2| $, it follows that $\sigma_2^\pm$ cut $T^\vee_2$ into two cylinders, one contains $\eta_1$, the other contains $\eta_2$. Consequently, $\sigma_1^\pm$ and $\sigma_2^\pm$ give rise to two pairs of homologous saddle connections in $\Sig$ which determine a special splitting $X'=(T'_1,T'_2,T'_3,v'_1,v'_2)$. We have

\begin{eqnarray*}
\Aa(T'_1) & = &\Aa(T^\vee_3)-|v'_1\wedge w_2|= \frac{A_2}{\alpha+1}+A_1-\frac{A_2}{\alpha+1}=A_1\\
\Aa(T'_2) & = & |(v'_1+v'_2)\wedge w_2| = A_2
\end{eqnarray*}

\noindent Therefore, $\Aa(T'_3)=A_3$. Since $\Lambda(T'_3)=\Z w_1 \oplus \Z(v'_2+\hat{u}_1)$, it follows from the choice of $\hat{u}_1$ that $v'_2$ is generic with respect to $\Lambda(T'_3)$.  We can then conclude that for any $(A_1,A_2,A_3,\alpha) \in \R^4_{>0}$ such that

\begin{itemize}
\item[.] $A_1+A_2+A_3=1$,

\item[.] $\DS{A_1+\frac{A_2}{\alpha+1}=A_1^0 +\frac{A_2^0}{\alpha_0+1}}$,
\end{itemize}

\noindent there exist $X'=(T'_1,T'_2,T'_3,v'_1,v'_2) \in \SSp(A_1,A_2,A_3,\alpha)$, with $v'_2$ generic with respect to $\Lambda(T'_3)$, such that $\Psi(X')\in \overline{\SL(2,\R)\cdot\Psi(X_0)}$. We can now complete the proof of Lemma \ref{PrfCLm2} as follows: first, for any $(A_1,A_2,A_3,\alpha)$ such that $\DS{\frac{A_1}{A_2}(\alpha+1)^2\notin \Q}$, it follows from Lemma \ref{PrfCLm1} that $\Psi(\SSp(A_1,A_2,A_3,\alpha)) \subset \overline{\SL(2,\R)\cdot\Psi(X_0)}$. By continuity of $\Psi$, it follows that $\overline{\SL(2,\R)\cdot\Psi(X_0)}$ contains $\Psi(\SSp(A_1,A_2,A_3,\alpha))$ for all $(A_1,A_2,A_3,\alpha)$. \carre

To complete the proof of Proposition \ref{PropC}, we need the following

\begin{lemma}\label{PrfCLm3}
For any $(A_1,A_2,A_3,\alpha)$ such that 

\begin{itemize}

\item[.] $A_1+A_2+A_3=1$,

\item[.] $\DS{A_1+\frac{A_2}{\alpha+1}< 1-(A_1^0+\frac{A_2^0}{\alpha_0+1})=\frac{\alpha_0A_2^0}{\alpha_0+1}+A_3^0}$,

\end{itemize}

\noindent we have $\Psi(\SSp(A_1,A_2,A_3,\alpha)) \subset \overline{\SL(2,\R)\cdot\Psi(X_0)}$.
\end{lemma}

\dem Since $\DS{A_1+\frac{A_2}{\alpha+1}< 1-(A_1^0+\frac{A_2^0}{\alpha_0+1})}$, we can find  $(A'_1,A'_2,A'_3,\alpha')\in \R_{>0}^4$ such that

\begin{itemize}
\item[.] $A'_1+A'_2+A'_3=1$,

\item[.] $A'_1+\frac{A'_2}{\alpha'+1}=A_1^0+\frac{A_2^0}{\alpha_0+1}$, and

\item[.]  $A'_3=A_1+\frac{A_2}{\alpha+1}$.

\end{itemize}

\noindent From Lemma \ref{PrfCLm2}, we know that  $\Psi(\SSp(A'_1,A'_2,A'_3,\alpha'))\subset \overline{\SL(2,\R)\cdot\Psi(X_0)}$. Consider an element $X=(T_1,T_2,T_3,v_1,v_2) \in \SSp(A'_1,A'_2,A'_3,\alpha')$, where  $v_1=(1,0), \; v_2=(\alpha',0)$.  Let $\Lambda_i$ denote the lattice associated to $T_i, i=1,2,3$.  Observe that we can choose $X$ such that (see Figure 6)

\begin{itemize}
\item[.] $\Lambda_1$ contains no vertical vectors,

\item[.] $\Lambda_2$ contains a vector vertical vector $w_2$ such that $\Lambda_2=\Z (v_1+v_2)\oplus \Z w_2$,

\item[.] $\Lambda_3=\Z v_3\oplus \Z w_1$, where $v_3$ is horizontal, and $w_1$ is vertical.
\end{itemize}

\begin{figure}[!ht]\label{Fig6}
\begin{center}
\begin{tikzpicture}[scale=0.5]

\draw (0,0) -- (0,2) -- (-4,2) -- (-4,0) -- (0,0) -- (7,0) -- (7,2) -- (11,2) -- (11,5) -- (12,7) -- (8,7) -- (7,5) -- (0,5) -- (0,2)-- (7,2) -- (7,5) -- (11,5);

\draw[>=angle 45, <->] (0,-1) -- (7,-1); \draw (3.5,-1) node[above] {$\alpha'$};

\draw[>=angle 45, <->] (7,-1) -- (11,-1); \draw (9,-1) node[above] {$1$};

\draw[>=angle 45, <->] (-4,-1.5) -- (7,-1.5); \draw (1.5,-1.5) node[below] {$\ell_3$};

\draw[>=angle 45, <->] (-4.5,0) -- (-4.5,2); \draw (-4.5,1) node[left] {$h_3$}; 

\draw[>=angle 45, ->] (12,0) -- (12,2); \draw (12,0) +(-0.2,0) -- +(0.2,0); \draw (12,1) node[right] {$w_1$};

\draw[>=angle 45, ->] (12,2) -- (12,5); \draw (12,2) +(-0.2,0) -- +(0.2,0); \draw (12,3.5) node[right] {$w_2$};

\draw (9,5.5) node[above] {$C_1$} (3.5,3) node[above] {$C_2$} (3.5,1) node {$C_3$};

\foreach \x in {(0,0),(0,2),(-4,2),(-4,0),(7,0),(7,2),(11,2),(11,5),(12,7),(8,7),(7,5),(0,5)} \filldraw[red] \x circle (2pt);

\end{tikzpicture}

\caption{Switching between horizontal and vertical splittings}
\end{center}
\end{figure}
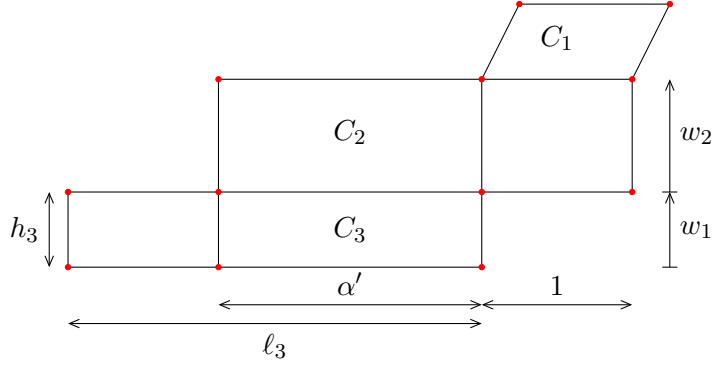

\noindent By assumption, we see that all $\Lambda_1, \Lambda_2, \Lambda_3$ contain horizontal vectors. Let $C_i, i=1,2,3,$ denote the horizontal cylinder obtained by slitting $T_i$ along the horizontal saddle connections, which correspond to the primitive horizontal vectors in $T_i$ . Let  $\ell_i$ and $h_i$ denote  width and the height of $C_i$. Note that $h_1,h_2$ are determined by $(A'_1,A'_2,A'_3,\alpha')$, and $\ell_3$ and $h_3$ must satisfy $\ell_3  >  \alpha' \text{ and } \ell_3 h_3 =  A'_3$.\\
By construction, the surface $\Sig$ constructed from $X$ admits another special splitting $X^\vee=(T^\vee_1,T^\vee_2,T^\vee_3,w_1,w_2)$ which is dual to $X$. Since $\Lambda_1$ contains no vertical vectors, the lattice $\Lambda^\vee_3$ does not contain any vertical vector. Let $\bar{m}^\vee$ denote the moduli ratio of $X^\vee$.

$$ \bar{m}^\vee=\frac{(\ell_3-\alpha')(h_2+h_3)}{\alpha'h_2}.$$

\noindent Since $\ell_3=A'_3/h_3$, we see that $\bar{m}^\vee$ is a non-constant rational function of $h_3$. Therefore, we can find $h_3$ so that $\bar{m}^\vee \notin \Q$. We deduce that there exists an element $X$ in $\SSp(A'_1,A'_2,A'_3,\alpha')$ such that the element $X^\vee$ defined above satisfies the conditions of Theorem \ref{ThC}.  Let $A_i^\vee, i=1,2,3,$ denote the area of $T^\vee_i$, and $\alpha^\vee=|w_2|/|w_1|$. By construction, we have 

$$A^\vee_1+\frac{A^\vee_2}{\alpha^\vee+1}=A'_3=A_1+\frac{A_2}{\alpha+1}.$$

\noindent Therefore, it follows from Lemma \ref{PrfCLm2} that 

$$\Psi(\SSp(A_1,A_2,A_3,\alpha)) \subset \overline{\SL(2,R)\cdot\Psi(X^\vee)}\subset \overline{\SL(2,\R)\cdot\Psi(X_0)}.$$ 

\carre

\subsection{Proof of Proposition \ref{PropC}}

All we need to show is that $\Psi(\SSp(A_1,A_2,A_3,\alpha)) \subset \overline{\SL(2,\R)\cdot \Psi(X_0)}$ for all $(A_1,A_2,A_3,\alpha)$ such that $A_1+A_2+A_3=1$. Choose $(A'_1,A'_2,A'_3,\alpha')$ in $\R_{>0}^4$ so that 

$$A'_1+\frac{A'_2}{\alpha'+1}< \min\{1-(A_1^0+\frac{A_2^0}{\alpha_0+1}), 1-(A_1+\frac{A_2}{\alpha+1})\}.$$

\noindent by Lemma \ref{PrfCLm3}, we know that $\Psi(\SSp(A'_1,A'_2,A'_3,\alpha')) \subset \overline{\SL(2,\R)\cdot \Psi(X_0)}$. Let $X$ be an element in $\SSp(A'_1,A'_2,A'_3,\alpha')$ which satisfies the conditions of Theorem \ref{ThC}. Since we have 

$$ A_1+\frac{A_2}{\alpha+1} < 1-(A'_1+\frac{A'_2}{\alpha'+1}),$$

\noindent by applying Lemma \ref{PrfCLm3} with $X$ in the place of $X_0$, we see that

$$\Psi(\SSp(A_1,A_2,A_3,\alpha)) \subset \overline{\SL(2,\R)\cdot \Psi(X)} \subset \overline{\SL(2,\R)\cdot\Psi(X_0)},$$

\noindent and the  proposition follows. \carre

\section{ Proof of Theorem \ref{ThC}}

By Proposition \ref{PropC}, we know that $\overline{\mathcal{O}}=\overline{\SL(2,\R)\cdot\Psi(X_0)}$ contains all the surfaces that admit a special splitting. We will show that $\overline{\mathcal{O}}$ contains all the Veech surfaces in $\H_1^\mathrm{hyp}(4)$, in particular, $\overline{\mathcal{O}}$ contains all the square-tiled surfaces. Since the set of square-tiled surfaces is dense in $\H^\mathrm{hyp}_1(4)$, it follows immediately that $\overline{\mathcal{O}}=\H^\mathrm{hyp}_1(4)$.\\
Let $\Sig$ be a Veech surface in $\H^\mathrm{hyp}_1(4)$. From Corollary \ref{SepSC}, we know that there exists on $\Sig$ a pair of homologous saddle connections $\delta^\pm$ such that by cutting along $\delta^\pm$, and gluing the two geodesic segments in the boundary of each of the connected component obtained from the cutting, we get a surface in $\H(0,0)$, which will be  denoted by $\Sig'$, and a surface in $\H(2)$ which will be denoted by $\Sig''$. On both $\Sig'$ and $\Sig''$ we have a marked saddle connection corresponding to the pair $\delta^\pm$, we denote both of them by $\delta$. Without loss of generality, we can assume that $\delta$ is horizontal. Since $\Sig$ is a Veech surface, $\Sig$ is decomposed into cylinders which are filled with horizontal closed geodesics. In particular, we see that $\Sig''$ is a union of horizontal cylinders. We have to possibilities

\begin{itemize}

\item[$\bullet$] \underline{\bf Case 1:} $\Sig''$ is the union of two cylinders. In this case, there exists another pair of homologous horizontal saddle connections $\gamma^\pm$ in $\Sig''$ which, together with $\delta^\pm$, determine a special splitting of $\Sig$. Therefore, $\Sig \in \overline{\mathcal{O}}$  by Proposition \ref{PropC}.

\item[$\bullet$] \underline{\bf Case 2:} $\Sig''$ contains only one horizontal cylinder. In this case, there exist two other horizontal saddle connections $\gamma_1, \gamma_2$ in $\Sig$ such that $\delta*\gamma_1*\gamma_2$ is  freely homotopic to a simple closed geodesic. Consequently, $\Sig''$ can be constructed from a single parallelogram $\P$  by the gluing as shown in  Figure 7. Actually, $\P$ is an octagon whose opposite sides are parallel and have the same length.  Let $U=(x,0),V_1=(y,0),V_2=(z,0)$, with $x >0, y>0, z>0,$ be the vectors associated to the saddle connections $\delta, \gamma_1,\gamma_2$ respectively. 

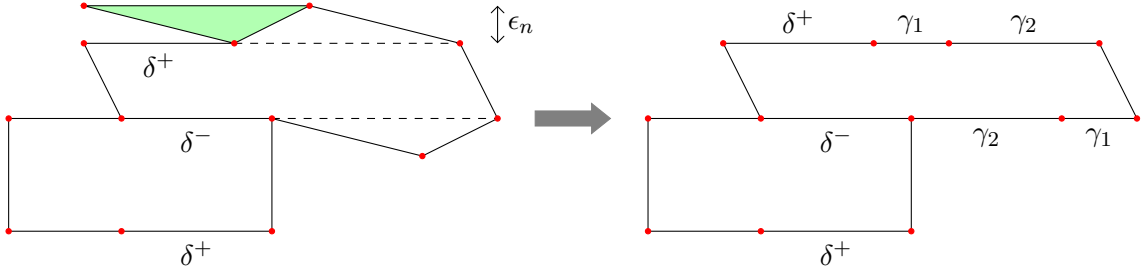
\begin{figure}[!ht]\label{Fig7}
\begin{center}
\begin{tikzpicture}[scale=0.5]

\fill[green!30] (-8,5) -- (-6,6) -- (-12,6) -- cycle;

\draw (-7,3) -- (-11,3) -- (-14,3) -- (-14,0) -- (-11,0) -- (-7,0) -- (-7,3) -- (-3,2) -- (-1,3) -- (-2,5) -- (-6,6) -- (-12,6) -- (-8,5) -- (-12,5) -- (-11,3) (-6,6) -- (-8,5);

\draw (10,3) -- (3,3) -- (3,0) -- (10,0) -- (10,3) -- (16,3) -- (15,5) -- (5,5) -- (6,3);

\filldraw[gray] (0,3) +(0,0.2) -- +(1.5,0.2) -- +(1.5,0.4) -- +(2,0) -- +(1.5,-0.4) -- +(1.5,-0.2) -- +(0,-0.2) -- cycle;

\draw[dashed] (-8,5) -- (-2,5) (-7,3) -- (-1,3);

\foreach \x in {(-7,3),(-11,3),(-14,3),(-14,0),(-11,0),(-7,0),(-3,2),(-1,3),(-2,5),(-6,6),(-12,6),(-8,5),(-12,5)} \filldraw[red] \x circle (2pt);
\foreach \x in {(3,0), (6,0), (10,0), (10,3), (14,3), (16,3), (15,5), (11,5), (9,5), (5,5), (6,3), (3,3)} \filldraw[red] \x circle (2pt);

\draw (-10,5) node[below] {$\delta^+$} (-9,3) node[below] {$\delta^-$} (-9,0) node[below]  {$\delta^+$}; 
\draw (7,5) node[above] {$\delta^+$} (8,3) node[below] {$\delta^-$} (8,0) node[below] {$\delta^+$};

\draw (10,5) node[above] {$\gamma_1$} (13,5) node[above] {$\gamma_2$} (15,3) node[below] {$\gamma_1$} (12,3) node[below] {$\gamma_2$};

\draw[>= angle 90, <->] (-1,5) -- (-1,6); \draw (-1,5.5) node[right] {$\epsilon_n$};

\end{tikzpicture}
\caption{Surfaces with special splitting converging to $\Sig$}
\end{center}
\end{figure}

\noindent Let $\{\epsilon_n\}$ be a sequence of positive real number decreasing to zero. For each $\epsilon_n$, we construct a surface $\Sig_n$ in $\H^\mathrm{hyp}_1(4)$ as follows: first, we construct a surface $\Sig''_n$ from an octagon $\P_n$, which is obtained from $\P$ by replacing $V_1$ by the vector $V^{(n)}_1=(y,\epsilon_n)$, and $V_2$ by the vector $V_2^{(n)}=(z,-\epsilon_n)$, then we glue $\Sig''_n$ to $\Sig'$ along the pair of homologous saddle connections $\delta^\pm$, and rescale to get a surface in $\H^\mathrm{hyp}_1(4)$. By construction, we see that $\Sig_n$ admits a special splitting by horizontal saddle connections, therefore $\Sig_n \in \overline{\mathcal{O}}$. As $\epsilon_n \lra 0$, the sequence $\{\Sig_n\}$ converges to $\Sig$, hence we have $\Sig \in \overline{\mathcal{O}}$.\\

\end{itemize}

The proof of Theorem \ref{ThC} is now complete. \carre

\section{Proof of Theorem \ref{ThB}}

We can now prove Theorem \ref{ThB} as a consequence of Theorem \ref{ThC}. The idea is to show that there exists in the closure of $\SL(2,\R)\cdot\Psi(X_0)$ a surface which admits a special splitting satisfying the conditions of Theorem \ref{ThC}.  As usual, let $A_i^0$ and $\Lambda_i^0$ denote the area and the associated lattice of $T_i^0, i=1,2,3$. We can assume that $v_1^0=(1,0)$ and $v_2^0=(\alpha_0,0)$. Let $t_1^0$ and $t_2^0$ denote the twists of the pairs $(T_1^0,v_1^0)$ and $(T_2^0,v_1^0+v_2^0)$ respectively (see \ref{SpSplSect}). Obviously, we only have to consider the case $\bar{m}_0=m_1^0/m_2^0 = (\alpha_0+1)^2A_1^0/A_2^0\in \Q$.\\

\noindent  Let $n_1, n_2$ be the integers such that $\gcd(n_1,n_2)=1$ and $n_1m_1^0+n_2m_2^0=0$. Applying Ratner's Theorem for the action of $U=\{\left(%
\begin{array}{cc}
  1 & t \\
  0 & 1 \\
\end{array}%
\right), \; \; t \in \R\}$, we see that $\overline{U\cdot\Psi(X_0)}$ contains $\Psi(X)$ for all $X=\SpEmt \in \SSp(A_1^0,A_2^0,A_3^0,\alpha_0)_\mathrm{hor}$ such that 

\begin{equation*}
n_1(t_1-t_1^0)+n_2(t_2-t_2^0) \in \Z
\end{equation*}

\noindent where $t_1, t_2$ are the twists of the pairs $(T_1,v_1)$ and $(T_2,v_1+v_2)$ respectively. Consider such an $X$ with $t_2=0$. Let $\Lambda_i$ denote the lattice associated to $T_i, i=1,2,3$. Since $t_2=0$, the lattice $\Lambda_2$ contains vertical vectors, let $w_1$ be the  primitive vertical vector in $\Lambda_2$, and let $\eta_2^\pm$ denote the pair of homologous saddle connections in $T_2$ such that $V(\eta_2^\pm)=w_2$. Let $w_1$ be a primitive vector in the lattice $\Lambda_3$ such that $|w_1\wedge v_2^0| < A_3^0$, and let $\eta_1^\pm$ denote the pair of homologous saddle connections in $T_3$ such that $V(\eta_1^\pm)=w_1$. The saddle connections $\eta_1^\pm$ and $\eta_2^\pm$ determine a splitting $X^\vee=(T_1^\vee,T_2^\vee, T_3^\vee, w_1,w_2)$ of the surface $\Sig= \Psi(X)$. Let $\Lambda_i^\vee$ and $A^\vee_i$ denote the associated lattice and the area of $T^\vee_i, i=1,2,3$. Here we have two cases:\\

\noindent \underline{\bf Case 1: $t_1\notin\Q$}\\
\noindent In this case, the lattice $\Lambda^\vee_3$ does not contain any vertical vector. We can choose  $w_1$ to be vertical and $\DS{\frac{A^\vee_1}{A^\vee_2}(\alpha^\vee+1)^2 \notin \Q}$, where $\alpha^\vee=|w_1|/|w_2|$ (see Lemma \ref{PrfCLm3}). Hence the splitting $X^\vee$ satisfies the condition of Theorem \ref{ThC}, it follows immediately that $\OrbClo=\overline{\SL(2,\R)\cdot\Psi(X^\vee)}=\H^\mathrm{hyp}_1(4)$.\\

\noindent \underline{\bf Case 2: $t_1 \in \Q$}\\
\noindent In this case, $\Lambda^\vee_3$ contains vertical vectors. Let $\hat{u}_3=(0,h_3)$, with $h_3>0$, be the primitive vertical  vector of $\Lambda^\vee_3$, and $\hat{v}_3=(\ell,h)$, with $\ell>0$ and $0\leq h<h_3$, be another primitive vector such that $\Lambda^\vee_3=\Z \hat{u}_3\oplus\Z \hat{v}_3$.\\
By assumption, we have $w_2=(0,h_2)$, with $h_2>0$. Remark that we have $h_3>h_2$. Recall that we are free to choose $T_3$ and $w_1$ provided $\Aa(T_3)=A_3^0$, and $|w_1\wedge v_2^0|<A_3^0$. By construction, we have $\Lambda_2^\vee=\Z v_2^0\oplus \Z(w_1+w_2)$.  The theorem follows from the following observation (see Lemma \ref{PrfCLm2})\\

\noindent {\bf Claim:} \textit{ We can choose $\Lambda_3$ and $w_1$ so that there exist a primitive vector $v'_1$ of $\Lambda^\vee_3$, and a primitive vector $v'_2$ of $\Lambda^\vee_2$ such that the surface $\Sig$ admits a special splitting $X'=(T'_1,T'_2,T'_3,v'_1,v'_2)$ dual to $X^\vee$ which satisfies the conditions of Theorem \ref{ThC}}.\\

\noindent \underline{\bf Proof of the claim:} Set $A_1=\ell(h_3-h_2)>0$, and choose $(A_2,A_3,\alpha)$ in $\R_{>0}^3$ such that

\begin{itemize}

\item[.] $A_1+A_2+A_3=1$,

\item[.] $\DS{A_1+\frac{A_2}{\alpha+1}=A_1^0+\frac{A_2^0}{\alpha_0+1}}$,

\item[.] $\DS{\frac{A_1}{A_2}(\alpha+1)^2 \notin \Q}$.

\end{itemize} 

%

\begin{figure}[!ht]\label{Fig8}
\begin{center}
\begin{tikzpicture}[scale=0.4]

\fill[gray!30] (7,15) -- (7,9) -- (intersection of 7,9 -- 11,21 and 7,15 -- 12,21) -- cycle;
\fill[gray!30] (12,15) -- (12,21) -- (7,6) -- (intersection of 7,6 -- 12,12 and 12,15 -- 10,9) -- cycle;

\fill[green!30] (7,6) -- (5,0) -- (8,0) -- cycle; \fill[green!30] (7,9) -- (4,9) -- (1,0) -- (4,0) -- cycle; \fill[green!30] (3,9) -- (0,9) -- (-1,6) -- (-4,5) -- (-3,-1) -- (0,0) -- cycle;

\draw (0,0) -- (-1,6) -- (-4,5) -- (-3,-1) -- (0,0) -- (8,0) -- (7,6) -- (12,12) -- (12,21) -- (7,15) -- (7,9) -- (-1,9) -- (-1,6)  (7,6) -- (7,9);

\draw[red] (7,9) -- (intersection of 7,9 -- 11,21 and 7,15 -- 12,21) (12,15) -- (intersection of 12,15 -- 10,9 and 7,6 -- 12,12);
\draw[red] (7,6) -- (12,21) (-1,6) -- (0,9) (0,0) -- (3,9) (1,0) -- (4,9)  (4,0) -- (7,9) (5,0) -- (7,6);

\draw[red] (8.5,15) node {$\sigma_1^+$} (10.2,16) node[right] {$\sigma_1^-$} (10.5,12.5) node {$\sigma_1^+$};
\draw[red] (5.5,6) node {$\sigma_2^+$} (6.3,2) node {$\sigma_2^-$} (1.5,6) node {$\sigma_2^+$} (2.3,2) node {$\sigma_2^-$};

\foreach \x in {(0,0), (-1,6), (-4,5), (-3,-1), (8,0), (7,6), (12,12), (12,15), (12,21), (7,15), (7,9), (-1,9)}
\filldraw[red] \x circle (2pt);

\draw[>= angle 45, <->, very thin] (7,6) -- (12,6); \draw (9.5,6) node[below] {$\ell$}; 
\draw[>= angle 45, <->, very thin] (12.5,6) -- (12.5,12); \draw (12.5,9) node[right] {$h$};
\draw[>= angle 45, <->, very thin] (12.5,12) -- (12.5,21); \draw (12.5,16.5) node[right] {$h_3$};
\draw[>= angle 45, <->, very thin] (-2,6) -- (-2,9); \draw (-2,7.5) node[left] {$h_2$};

\draw[>= angle 45, ->] (-4,-1) -- (-5,5); \draw (-4.5, 2) node[left] {$w_1$}; \draw (-4,-1) +(0.3,0.05) -- +(-0.3,-0.05);
\draw[>= angle 45, ->] (-5,6) -- (-5,9); \draw (-5,7.5) node[left] {$w_2$}; \draw (-5,6) +(-0.2,0) -- +(0.2,0);

\draw[>=angle 45, ->] (-3,-2) -- (0,-1); \draw (-3,-2) +(0.05,-0.15) -- +(-0.05,0.15); \draw (-1.5,-1.5) node[below] {$\hat{u}_1$};

\draw[>=angle 45, ->] (-1,9.5) -- (7,9.5); \draw (-1,9.5) +(0,0.2) -- +(0,-0.2); \draw (4,9.5) node[above] {$v_2^0$};

\end{tikzpicture}

\caption{{\bf Case 2} the lattice $\Lambda^\vee_3$ contains vertical vectors.}

\end{center}
\end{figure}

\noindent Since $\DS{A_1+\frac{A_2}{\alpha+1}=A_1^0+\frac{A_2^0}{\alpha_0+1}}=A^\vee_3=\ell h_3$, and $\DS{A_1=\ell(h_3-h_2)}$, it follows $\DS{\frac{A_2}{\alpha+1}=\ell h_2}$, hence

\begin{equation}\label{PrfBIneq1}
\frac{\alpha A_2}{\alpha+1}=\alpha\ell h_2<1-(A_1^0+\frac{A_2^0}{\alpha_0+1})=\frac{\alpha_0A_2^0}{\alpha_0+1}+A_3^0=A_0
\end{equation}

\noindent Choose $q \in \N$ large enough so that

$$\left\{%
\begin{array}{l}
    \DS{\frac{h}{h_3q}< \frac{\alpha\ell h_2}{\alpha\alpha_0 h_3},} \\
    
    \DS{\frac{1}{q} < \frac{1}{2}\frac{A_0-\alpha\ell h_2}{\alpha\alpha_0h_3}}.
\end{array}%
\right.$$

\noindent From (\ref{PrfBIneq1}), it follows  that there exists $p\in \N$ such that

$$ \frac{\alpha\ell h_2}{\alpha\alpha_0h_3}-\frac{h}{h_3q} < \frac{p}{q} < \frac{A_0}{\alpha\alpha_0h_3}-\frac{h}{h_3q}.$$

\noindent Now, we can take

\begin{itemize}
\item[.] $v'_1=p\hat{u}_3+\hat{v}_3=(\ell, ph_3+h)$,

\item[.] $\DS{w_1=(\frac{\alpha \ell-\alpha_0}{q}, \alpha h_3(\frac{h}{h_3q}+\frac{p}{q})-h_2)}$,

\item[.] $v'_2=v_2^0+q(w_1+w_2)=(\alpha\ell, \alpha(ph_3+h))$.

\end{itemize}

\noindent Observe that $v'_1$ and $v'_2$ are primitive vectors in $\Lambda^\vee_3$ and $\Lambda^\vee_2$ respectively. Clearly, we have $v'_2=\alpha v'_1$.  By the choice of $p,q$, we also have

\begin{itemize}
\item[.] $\DS{|v'_1\wedge w_2|=\ell h_2 < \ell h_3=A^\vee_3=A_1^0+\frac{A_2^0}{\alpha_0+1}}$,

\item[.] $\DS{|v'_2\wedge w_2|=\alpha\ell h_2 < \alpha\alpha_0h_3(\frac{h}{h_3q}+\frac{p}{q})=|v^0_2\wedge(w_1+w_2)|= A^\vee_2}$,

\item[.] $\DS{A^\vee_2=\alpha\alpha_0h_3(\frac{h}{h_3q}+\frac{p}{q})<A_0=1-A^\vee_3}$,



\end{itemize}

\noindent Consequently, the surface $\Sig$ admits a special splitting determined by two pairs of homologous saddle connections $\sigma_1^\pm$ and $\sigma_2^\pm$, where $\sigma_1^\pm$ is the pair of saddle connections in $T^\vee_3$ corresponding to $v'_1$, and $\sigma_2^\pm$ is the pair of saddle connections in $T^\vee_2$ corresponding to $v'_2$ (see Figure 8). Let $X'=(T'_1,T'_2,T'_3,v'_1,v'_2)$ denote this special splitting, then we have $\Aa(T'_i)=A_i, i=1,2,3$.\\
By construction, $w_1$ is a primitive vector of $\Lambda^\vee_1$. Let $\hat{u}_1$ be another primitive such that $\Lambda^\vee_1=\Z w_1\oplus \Z \hat{u}_1$. Recall that we can choose $T_3$ arbitrarily provided $\Aa(T_3)=A_3^0$, therefore we are free to choose $\hat{u}_1$, provided $|\hat{u}_1\wedge w_1|=A^\vee_1=1-(A^\vee_2+A^\vee_3)$. It is easy to check that we can choose such a $\hat{u}_1$ so that the vector $v'_2$ is generic with respect to $\Lambda(T'_3)=\Z w_1\oplus\Z(v'_2+\hat{u}_1)$. With this choice, we see that splitting $X'$ satisfies the conditions of Theorem \ref{ThC}, and the claim is then proved. \carre

By Theorem \ref{ThC}, we know that $\overline{\SL(2,\R)\cdot\Sig}=\H^\mathrm{hyp}_1(4)$. Since $\Sig \in \overline{\SL(2,\R)\cdot\Sig_0}$, it follows $\overline{\SL(2,\R)\cdot\Sig_0}=\H^\mathrm{hyp}_1(4)$. The proof of Theorem \ref{ThB} is now complete. \carre

\section{Surfaces admitting completely periodic directions with three cylinders }

\subsection{Two models of decomposition into three cylinders}

\begin{lemma}\label{CylDecLm}

Let $\Sig$ be a surface in $\H^\mathrm{hyp}(4)$. Assume that $\Sig$ is decomposed into three horizontal cylinders, that is, the horizontal direction is completely periodic for $\Sig$ with three cylinders. Then the surface $\Sig$ can be reconstructed from three (horizontal) cylinders by one of the following gluing models

\begin{figure}[ht]\label{Fig9}
\begin{center}
\begin{tikzpicture}[scale=0.4]

\draw (-14,4) -- (-16,4) -- (-16,0) -- (-10,0) -- (-10,4) -- (-6,4) -- (-6,7) -- (-5,9) -- (-9,9) -- (-10,7) -- (-14,7) -- (-14,4) -- (-10,4) (-10,7) -- (-6,7);
\draw (16,5) -- (9,5)  (9,2) -- (3,2) -- (0,2) -- (0,5) -- (9,5) -- (8,8) -- (15,8) -- (16,5) -- (16,2) -- (9,2) -- (10,0) -- (4,0) -- (3,2);


\draw[dashed] (-14,0) -- (-14,4) (-10,4) -- (-10,7);

\draw (-7,9) +(0,-0.2) -- +(0,0.2); \draw (-8,4) +(0,-0.2) -- +(0,0.2); 
\draw (-12,7) +(-0.1,-0.2) -- +(-0.1,0.2) +(0.1,-0.2) -- +(0.1,0.2); \draw (-12,0) +(-0.1,-0.2) -- +(-0.1,0.2) +(0.1,-0.2) -- +(0.1,0.2);
\draw (-15,4) +(-0.2,-0.2) -- +(-0.2,0.2) +(0,-0.2) -- +(0,0.2) +(0.2,-0.2) -- +(0.2,0.2); \draw (-15,0) +(-0.2,-0.2) -- +(-0.2,0.2) +(0,-0.2) -- +(0,0.2) +(0.2,-0.2) -- +(0.2,0.2);

\draw (11.5,8) +(-0.2,-0.2) -- +(-0.2,0.2) +(0,-0.2) -- +(0,0.2) +(0.2,-0.2) -- +(0.2,0.2); \draw (12.5,2) +(-0.2,-0.2) -- +(-0.2,0.2) +(0,-0.2) -- +(0,0.2) +(0.2,-0.2) -- +(0.2,0.2);
\draw (3,5) +(-0.1,-0.2) -- +(-0.1,0.2) +(0.1,-0.2) -- +(0.1,0.2); \draw (7,0) +(-0.1,-0.2) -- +(-0.1,0.2) +(0.1,-0.2) -- +(0.1,0.2);
\draw (7.5,5) +(0,-0.2) -- +(0,0.2); \draw (1.5,2) +(0,-0.2) -- +(0,0.2); 

\foreach \x in {(-14,4), (-16,4), (-16,0), (-14,0), (-10,0), (-10,4), (-6,4), (-6,7), (-5,9), (-9,9), (-10,7), (-14,7)}
\filldraw[red] \x circle (2pt);

\foreach \x in {(9,5), (6,5), (0,5), (0,2), (3,2), (4,0), (10,0), (9,2), (16,2), (16,5), (15,8), (8,8)} \filldraw[red] \x circle (2pt);


\draw (-12,-2) node {\bf Case I};
\draw (7,-2) node {\bf Case II};

\end{tikzpicture}
\caption{Two models of gluing}
\end{center}
\end{figure}
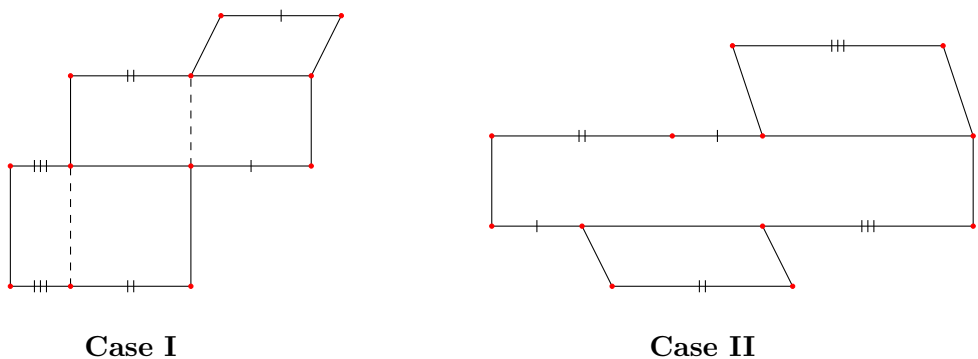

\end{lemma}

\dem First, observe that $\Sig$ has exactly $5$ horizontal saddle connections, since the angle at the unique singular point of $\Sig$ is $10\pi$. Let $C_1,C_2,C_3$ denote the three horizontal cylinders. Since each of the horizontal saddle connections is contained in the lower boundary component of a unique cylinder, we then have a partition of set of horizontal saddle connections into three subsets, there are only two such partitions corresponding to two ways of writing $5$ as the sum of three positive integers: $5=1+1+3=1+2+2$.\\

Next, let us show that the hyper-elliptic involution $\tau$ of $\Sig$ preserves each of the cylinders $C_i, i=1,2,3$. Consider a simple closed geodesic $c_i$ in $C_i$ close to its lower boundary. Since $\tau(c_i)+c_i=0$ in $H_1(\Sig,\Z)$, we deduce that $c_i$ and $\tau(c_i)$ cut $\Sig$ into two connected components, each of which is equipped with a flat metric with geodesic boundaries. Since $\Sig$ has only one singularity, one of the two components has no singularities in the interior, and must be a cylinder. Therefore, $c_i$ and $\tau(c_i)$ are contained in the same cylinder $C_i$. As a consequence, we see that $\tau$ maps the lower boundary of each cylinder to its upper boundary. In particular, the upper boundary and the lower boundary of each cylinder contain the same number of saddles connections, and moreover, for each saddle connection in the lower boundary is paired up with a saddle connection in the upper boundary, which is its image under $\tau$.\\

\noindent From these two observations, it is now easy to check that there are only two ways to construct $\Sig$ from three cylinders, which are shown in Figure 9. \carre 

\rem The fact that the hyper-elliptic involution preserves each of the cylinders is already known to Kontsevich-Zorich (see \cite{KonZo}, Lemma 8).

\subsection{Proof of Corollary \ref{CorB}}

\subsubsection{Proof of Corollary \ref{CorB}, Case I)}\label{PrfCorB1}

In this case, let $C_1^0$ denote the unique simple cylinder of the decomposition, $C_2^0$ denote the cylinder adjacent to $C_1^0$, and $C_3^0$ the remaining cylinder. Let $\ell_i^0,h_i^0, m_i^0$ denote respectively the width, the height, and the modulus of $C^0_i$.  Since $m_1^0,m_2^0,m_3^0$ are independent over $\Q$, by applying Ratner's Theorem for the action of $U$, we deduce that $\overline{U\cdot\Sig_0}$ contains all the surfaces $\Sig$ which are constructed from $3$ horizontal cylinders $C_1,C_2,C_3$ by the same gluing model, whenever $C_i$ has the same width and height as $C_i^0$.\\
On each boundary component of $C_1$ we have a marked point which corresponds to the unique singularity of $\Sig$. Let $v_1=(v_1^x,v_1^y)$ be the associated vector of any geodesic segment joining the marked point in the lower boundary to the marked point in the upper boundary. We then define the twist $t_1$ of $C_1$ to be $\DS{\frac{v_1^x}{\ell_1^0} \mod \Z}$. On each boundary component  of $C_2$ (resp. $C_3$), we have two marked points, therefore each boundary component is the union of two geodesic segments. From Lemma \ref{CylDecLm}, we see that each segment in the upper boundary of $C_2$ is paired up with a segment in the lower boundary component by the hyper-elliptic involution. Take such a pair of segments, and consider  a segment joining the left endpoint of the segment in the lower boundary to the left endpoint of the segment in the upper boundary. Let $v_2=(v_2^x,v_2^y)$ be the vector associated to this segment, we then define the twist $t_2$ of $C_2$ to be $\DS{t_2=\frac{v_2^x}{\ell_2^0}\mod \Z}$. We define the twist $t_3$ of $C_3$ in the same manner.\\
Observe that any value of $(t_1,t_2,t_3)$ gives us a unique surface $\Sig$ in $\overline{U\cdot\Sig_0}$. Consider the case $t_2=t_3=0$, in that case $\Sig$ admits a special splitting by two pairs of vertical homologous saddle connections. It is easy to see that if $t_1$ is not in $\Q$ then this splitting satisfies the condition of Theorem \ref{ThB}, that is the lattice associated to $T_3$ does not contain any vertical vector. It follows immediately that $\overline{\SL(2,\R)\cdot\Sig_0}=\overline{\SL(2,\R)\cdot\Sig}=\H^\mathrm{hyp}_1(4)$. 

\subsubsection{ Proof of Corollary \ref{CorB} Case II)}

In this case, we have two simple cylinders, which will be denoted by $C_1^0$ and $C_2^0$, the remaining cylinder has $3$ saddle connections in each boundary component, and will be denoted by  $C_3^0$. Let $\gamma_1^+,\gamma_2^+,\gamma_3^+$ denote the saddle connections contained in the upper boundary of $C_3^0$ such that $\gamma_1^+$ (resp. $\gamma_2^+$) is also the lower boundary of $C_1^0$ (resp. $C_2^0$). Let $\gamma_i^-$ denote the  image of $\gamma_i^+$ under $\tau$. Note that the lower boundary of $C^0_3$ is the union of $\gamma_1^-,\gamma_2^-,\gamma_3^-$, and in fact $\gamma_3^+=\gamma_3^-$ (see Figure 10). \\
Let $\ell_i, h_i$ denote respectively the width and the height of $C_i^0, i=1,2,3$. Since the cylinders $C_1^0$ and $C_2^0$ are simple, we define their twists $t_1^0,t_2^0$ as in Case I). Let $\delta^+$ (resp. $\delta^-$) denote a pair of homologous saddle connections in $C_3^0$ which joins the left (resp. right) endpoint of $\gamma_1^-$ to the left (resp. right) endpoint of $\gamma^+_1$. Using the action of $U$, we can assume that $\delta^\pm$ are vertical.\\
Applying the Ratner's Theorem, we see that $\overline{U\cdot\Sig_0}$ contains all surfaces obtained from three cylinders $(C_1,C_2,C_3)$ by the same gluing model as $(C_1^0,C^0_2,C_3^0)$, provided $C_i$ has the same width and height as $C_i^0$. In particular, $\overline{U\cdot\Sig_0}$ contains all surfaces constructed from three cylinders $(C_1,C_2,C_3)$ with $C_3=C_3^0$, and, for $i=1,2$, the twist $t_i$ of $C_i$ can be chosen arbitrarily. Let $\Sig$ be such a surface. Cut $\Sig$ along $\delta^\pm$, then glue the geodesic segments corresponding to $\delta^\pm$ on each component together, we get a surface in $\H(0,0)$ containing $C_1$, which will be denoted by $\Sig'$, and a surface in $\H(2)$ containing $C_2$, which will be denoted by $\Sig''$. In both of $\Sig'$ and $\Sig''$, we have a marked saddle connection corresponding to $\delta^\pm$, we denote both of them by $\delta$.

\begin{figure}[ht]\label{Fig10}
\begin{center}
\begin{tikzpicture}[scale=0.5]

\draw (16,5) -- (9,5) (9,2) -- (3,2) -- (0,2)  (0,5) -- (9,5) -- (8,8) -- (15,8) -- (16,5)  (16,2) -- (9,2) -- (10,0) -- (4,0) -- (3,2);
\draw[dashed, red] (0,2) -- (3,5) (4,0) -- (12,8) (7,0) -- (15,8) (13,2) -- (16,5);

\draw (0,5) -- (0,2) (9,5) -- (9,2) (16,5) -- (16,2);
\draw (0,3.5) node[left] {$\delta^+$} (9,3) node[left] {$\delta^-$} (16,3.5) node[right] {$\delta^+$};

\foreach \x in {(9,5), (6,5), (0,5), (0,2), (3,2), (4,0), (10,0), (9,2), (16,2), (16,5), (15,8), (8,8)} \filldraw[red] \x circle (2pt);

\draw (13.5,5) node[below] {$\gamma_1^+$} (7.5,5) node[above] {$\gamma_3^+$} (1.5,5) node[above] {$\gamma_2^+$};

\draw (1.5,2) node[below] {$\gamma_3^-$} (4.5,2) node[above] {$\gamma_2^-$} (11,2) node[below] {$\gamma_1^-$};

\draw[red] (10,6) node[above] {$\sigma_1^+$} (13,6) node[above] {$\sigma_1^-$} (15,4) node[below] {$\sigma_1^+$};

\draw[red] (7,3) node[above] {$\sigma_2^+$} (2,4) node[below] {$\sigma_2^-$};

\draw[very thin, >= angle 45, <->] (17.5,5) -- (17.5,8); \draw[very thin, >=angle 45, <->] (17.5,2) -- (17.5,5); \draw[very thin, >=angle 45, <->] (17.5,0) -- (17.5,2);  

\draw (17.5,6.5)  node[right] {$h_1$} (17.5, 3.5) node[right] {$h_3$} (17.5,1) node[right] {$h_2$};

\draw[very thin, >= angle 45, <->] (8,8.2) -- (15,8.2); \draw (11.5,8.2) node[above] {$\ell_1$};
\draw[very thin, >= angle 45, <->]  (4,-0.4) -- (10,-0.4); \draw (7,-0.4) node[below] {$\ell_2$};

\draw[very thin, >= angle 45, <->] (0,-0.2) -- (16,-0.2); \draw (12,-0.2) node[above] {$\ell_3$};

\end{tikzpicture}
\caption{Finding new special splittings}
\end{center}
\end{figure}
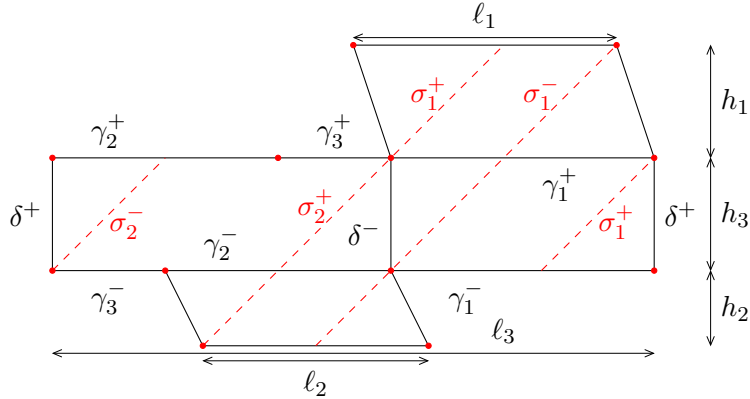

\noindent In $\Sig'$, for any $t_1 \in ]0,1[$, we have a pair of homologous saddle connections $\sigma_1^\pm$ which correspond to the vector $w_1=(t_1\ell_1,h_1+h_3)$. This pair of saddle connections cut $\Sig'$ into two cylinders, one of which contains $\delta$. \\
Suppose that $t_2\in ]0,1[$, then there exists a pair of homologous saddle connections $\sigma_2^\pm$ in $\Sig''$ $V(\sigma_2^\pm)= w_2=(t_2\ell_2, h_2+h_3)$, which bound a simple cylinder containing $\delta$. If we cut off the simple cylinder bounded by $\sigma_2^\pm$ from $\Sig''$, and then glue the geodesic segments corresponding to $\sigma^\pm$, we obtain a torus $T$ in $\H(0,0)$ together with a marked saddle connection $\sigma_2$. Let $\Lambda$ denote the lattice in $\R^2$ associated to $T$, then $\Lambda$ is generated by $u=(\ell_2,h_3)$ and $v=(t_2\ell_2-|\gamma_3^+|, h_2+h_3)$. Note that $u$ is independent of $t_2$.\\
Recall that $w_2$ is parallel to a vector in $\Lambda=\Z u\oplus\Z v$ if and only if we can write $v=\lambda u+ \mu w_2$ with $\lambda\in \Q$. As $t_2$ varies, we see that the set of $t_2$ for which $w_2$ is parallel to a vector in $\Lambda$ is countable, which means that, given any $\epsilon>0$, we can find $t_2 \in ]0,\epsilon[$ such that $w_2$ is not parallel to any vector in $\Lambda$. Therefore, we can find $\DS{t_2 \in ]0, \frac{\ell_1(h_1+h_3)}{\ell_2(h_2+h_3)}[}$, such that $w_2$ is not parallel to any vector in $\Lambda$. Now, take $\DS{t_1=\frac{\ell_2(h_2+h_3)}{\ell_1(h_1+h_3)}t_2}$, we have $t_1\in ]0;1[$, hence we can find $\sigma_1^\pm$ as above. By the choice of $t_1$ and $t_2$, $w_1$ and $w_2$ are parallel. Reconstruct $\Sig$ from $\Sig'$ and $\Sig''$, we see that $\sigma_1^\pm$ and $\sigma_2^\pm$ determine a special splitting of $\Sig$, which satisfies the condition of Theorem \ref{ThB}. Since $\Sig \in \overline{U\cdot\Sig_0}$, the corollary follows. \carre

\section{Applications}\label{ExSect}

\subsection{Generic surfaces with coordinates in a quadratic field}

In \cite{HubLanMo3}, Hubert-Lanneau-Möller raise the following question: does there exist a generic translation surface of genus $g$ with all coordinates in  a number field $K$ such that $[K:\Q]<g$? Theorem \ref{ThB} provides us with an affirmative answer to this question for the case $\H^\mathrm{hyp}(4)$. For every quadratic field $K$, one can construct a surface in $\H^\mathrm{hyp}_1(4)$ with all coordinates in $K$ which satisfies the condition of Theorem \ref{ThB}. Here below is such a surface with coordinates in $\Q[\sqrt{2}]$.

\begin{figure}[ht]\label{Fig11}
\begin{center}
\begin{tikzpicture}[scale=0.4]
\draw (6,8) -- (3,8) -- (3,4) -- (0,7) -- (0,3) -- (3,0) -- (6,0) -- (6,4) -- (11,4) -- (11,8) -- (11,10) -- (6,10) -- cycle; 

\draw[dashed] (6,8) -- (11,8) (3,4) -- (6,4);

\draw (6,9) +(-0.2,0) -- +(0.2,0) (11,9) +(-0.2,0) -- +(0.2,0);
\draw (3,6) +(-0.2,0.1) -- +(0.2,0.1) +(-0.2,-0.1) -- +(0.2,-0.1); \draw (11,6) +(-0.2,0.1) -- +(0.2,0.1) +(-0.2,-0.1) -- +(0.2,-0.1);
\draw (0,5) +(-0.2,0.2) -- +(0.2,0.2) +(-0.2,0) -- +(0.2,0) +(-0.2,-0.2) -- +(0.2,-0.2); \draw (6,2) +(-0.2,0.2) -- +(0.2,0.2) +(-0.2,0) -- +(0.2,0) +(-0.2,-0.2) -- +(0.2,-0.2); 

\draw (1.5,5.5) +(0,-0.2) -- +(0,0.2); \draw (1.5,1.5) +(0,-0.2) -- +(0,0.2); 
\draw (4.5, 8) +(-0.1,-0.2) -- +(-0.1,0.2) +(0.1,-0.2) -- +(0.1,0.2); \draw (4.5,0) +(-0.1,-0.2) -- +(-0.1,0.2) +(0.1,-0.2) -- +(0.1,0.2); 
\draw (8.5,10) +(-0.2,-0.2) -- +(-0.2,0.2) +(0,-0.2) -- +(0,0.2) +(0.2,-0.2) -- +(0.2,0.2); \draw (8.5,4) +(-0.2,-0.2) -- +(-0.2,0.2) +(0,-0.2) -- +(0,0.2) +(0.2,-0.2) -- +(0.2,0.2);

\foreach \x in {(6,8), (3,8), (3,4), (0,7), (0,3), (3,0), (6,0), (6,4), (11,4), (11,8), (11,10), (6,10)} \filldraw[red] \x circle (2pt);

\draw[>=angle 45, very thin, <->] (12,8) -- (12,10); \draw[>=angle 45, very thin, <->] (12,8) -- (12,4); \draw[>=angle 45, very thin, <->] (12,4) -- (12,0);
\draw[>=angle 45, very thin, <->] (0,-1) -- (3,-1); \draw[>=angle 45, very thin, <->] (3,-1) -- (6,-1); \draw[>=angle 45, very thin, <->] (6,-1) -- (11,-1);

\draw[>=angle 45, very thin, <->] (0,0) -- (0,3);

\draw (12,9) node[right] {$\frac{1}{4}$} (12,6) node[right] {$\frac{1}{2}$} (12,2) node[right] {$\frac{1}{2}$};
\draw (1.5,-1) node[below] {$\frac{\sqrt{2}}{4}$} (4.5,-1) node[below] {$\frac{2-\sqrt{2}}{2}$} (8.5,-1) node[below] {$\frac{\sqrt{2}}{2}$};
\draw (0,1.5) node[left] {$\frac{\sqrt{2}}{4}$};

\end{tikzpicture}

\caption{ A generic surface in $\H_1^\mathrm{hyp}(4)$ with coordinates in $\Q[\sqrt{2}]$}
\end{center}
\end{figure}
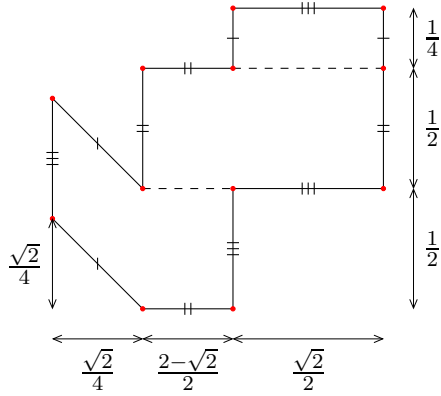

\subsection{Thurston-Veech surface with cubic trace field}

Surfaces obtained by the Thurston-Veech construction have large Veech group, which contains infinitely many hyperbolic elements (see \cite{HubLanMo3} for definition and further detail on Thurston-Veech construction). Recall that the trace field of a translation surface is the field generated over $\Q$ by the the traces of the matrices in its Veech group. If $K$ is the trace field of a translation surface of genus $g$ then $[K:\Q]\leq g$. For $g=2$, McMullen (\cite{McM2}) shows that if $[K:\Q]=2$ then the $\SL(2,\R)$-orbit of the surface can not be dense in its stratum. However, for $g=3$, Hubert-Lanneau-Möller (\cite{HubLanMo1}, \cite{HubLanMo2}) show that there exist surfaces in the hyper-elliptic locus $\mathcal{L}$ of $\H^\mathrm{odd}(2,2)$ obtained by the Thurston-Veech construction with trace field of degree $3$ such that the $\SL(2,\R)$-orbit is dense in $\mathcal{L}$. Note that $\mathcal{L}$ is a closed $\SL(2,\R)$-invariant subset of $\H^\mathrm{odd}(2,2)$, therefore these surfaces can be viewed as generic.\\

The surfaces obtained from Thurston-Veech construction are completely algebraically periodic in the sense of Calta-Smillie (see \cite{CalSm}). In particular, if such a surface admits a special splitting $\SpEmt$, then $v_2$ must be parallel to a vector in the lattice associated to $T_3$. Therefore, a Thurston-Veech surface can never satisfy the condition of Theorem \ref{ThB}. However, if the trace field is of degree $3$ over $\Q$, one can find examples of Thurston-Veech surfaces which admit decompositions into three parallel cylinders whose moduli are independent over $\Q$. By Corollary \ref{CorB}, it follows that the $\SL(2,\R)$-orbits of such surfaces are dense in $\H^\mathrm{hyp}_1(4)$. Here below, we will give the explicit construction of some of such surfaces.\\
We construct surfaces in $\H^\mathrm{hyp}(4)$ for which the horizontal and vertical directions are completely periodic with three cylinders. To get such a surface, we glue three horizontal cylinders $C_1,C_2,C_3$ as shown in Figure 12.\\

\begin{figure}[ht]\label{Fig12}
\begin{center}
\begin{tikzpicture}[scale=0.6]

\fill[blue!30] (10,8) -- (12,10) -- (9,10) -- cycle; \fill[blue!30] (10,4) -- (13,4) -- (17,8) -- (16,10) -- cycle;
\fill[green!30] (0,0) -- (2,0) -- (10,8) -- (6,8) -- (2,4) -- (0,4) -- cycle; \fill[green!30] (6,0) -- (10,0) -- (10,4) -- cycle;

\draw (0,0) -- (10,0) (0,4) -- (17,4) (2,8) -- (17,8) (9,10) -- (16,10) (9,10) -- (10,8) (16,10) -- (17,8) (2,8) -- (2,4) (10,8) -- (10,4) (17,8) -- (17,4) (0,4) -- (0,0) (2,4) -- (2,0) (10,4) -- (10,0);

\draw[red] (0,2) -- (6,8) (0,0) -- (8,8) (2,0) -- (12,10) (6,0) -- (16,10) (8,0) -- (10,2) (13,4) -- (17,8);

\draw (13,10) node[above] {$\delta_1^-$} (13,8) node[above] {$\delta_1^+$} (14.5,4) node[below] {$\delta_1^-$};
\draw (4.5, 8) node[above] {$\delta_2^+$} (8.5,4) node[above] {$\delta_2^-$} (4.5,0) node[below] {$\delta_2^+$};
\draw (1,4) node[above] {$\delta_3$} (1,0) node[below] {$\delta_3$};

\draw (0,3) node[left] {$\eta_1^+$} (10,3) node[right] {$\eta_1^+$};
\draw (2,6) node[left] {$\eta_2^-$} (10,6.5) node[right] {$\eta_2^+$} (17,6) node[right] {$\eta_2^-$}; 
\draw (9.5,9) node[left] {$\eta$} (16.5,9) node[right] {$\eta$}; 

\draw[red] (11,9) node[above] {$\sigma_1^+$} (12.5, 6.5) node[above] {$\sigma_1^-$} (15,6) node[below] {$\sigma_1^+$};
\draw[red] (5,3) node[below] {$\sigma_2^+$} (8,2) node[above] {$\sigma_2^-$} (9,1) node[above] {$\sigma_3$};
\draw[red] (5,5) node[above] {$\sigma_3$} (4,6) node[above] {$\sigma_2^-$};

\foreach \x in {(0,0), (2,0), (10,0), (10,4), (17,4), (17,8), (16,10), (9,10), (10,8), (2,8), (2,4), (0,4)} \filldraw[red] \x circle (2pt);

\draw[very thin, >=angle 45, <->] (19,0) -- (19,4); \draw[very thin, >=angle 45, <->] (19,4) -- (19,8); \draw[very thin, >=angle 45, <->] (19,8) -- (19,10);

\draw (19,2) node[right] {$h_3$} (19,6) node[right] {$1$} (19,9) node[right] {$h_1$};

\draw[very thin, >=angle 45, <->] (0,-2) -- (10,-2); \draw[very thin, >=angle 45, <->] (10,-2) -- (17,-2);
\draw[very thin, >=angle 45, <->] (2,-2.2) -- (17,-2.2);
\draw (7,-2) node[above] {$\ell_3$} (14,-2) node[above] {$\ell_1$} (9,-2.2) node[below] {$1$};

\end{tikzpicture}
\caption{Cubic Thurston-Veech surface with a non-parabolic completely periodic direction}
\end{center}
\end{figure}
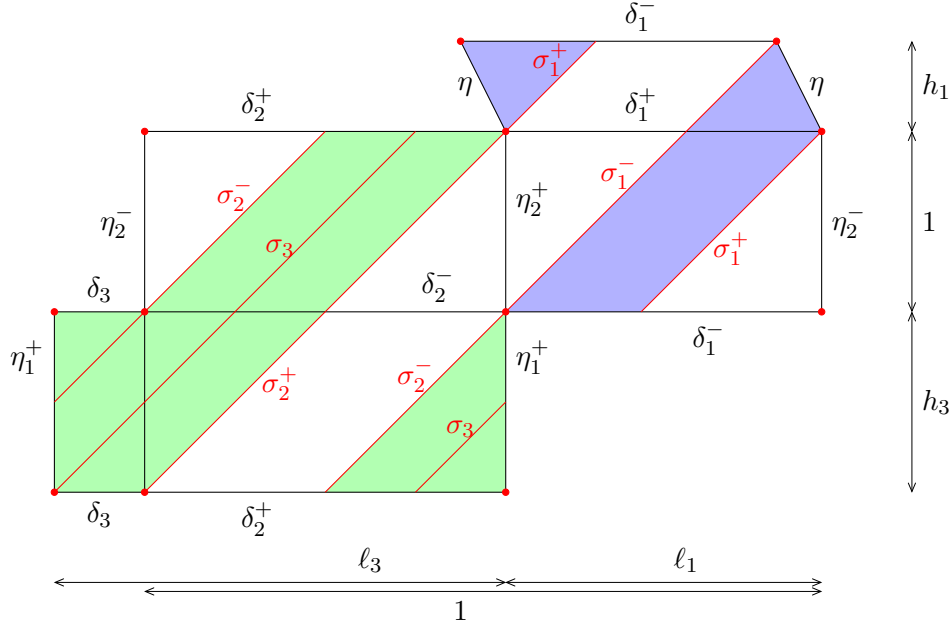

\noindent Let $\ell_i$ and $h_i$ denote the width and the height of $C_i, i=1,2,3$. We define the twists $t_i$ of $C_i, i=1,2,3,$ as in Section \ref{PrfCorB1}.  In what follows we fix $t_2=t_3=0$, and consider the cases $t_1=(n-1)/n, n\in \N$.\\
We denote by $\delta_1^\pm, \delta_2^\pm, \delta_3$ the  horizontal saddle connections contained in the boundary of $C_1,C_2,C_3$ as shown in Figure 12. We choose the orientation for every horizontal saddle connection to be from the left to the right, and for every vertical saddle connection to be upward. Let $\eta_2^+$ and $\eta_2^-$ (resp. $\eta_1^+$ and $\eta_1^-$) denote the vertical saddle connections in $C_2$ (resp. in $C_3$) joining the left and the right endpoints of $\delta_2^-$ to the left and the right endpoints of $\delta_2^+$ respectively. We see that the surface $\Sig$ admits a special spitting $(T^\vee_1, T^\vee_2, T^\vee_3, w_1,w_2)$ determined by $\eta^\pm_1$ and $\eta^\pm_2$ ($w_i=V(\eta^\pm_i), i=1,2$).  Since $t_1=(n-1)/n$, there exists a vertical saddle connection in $T^\vee_3$ which crosses $\delta_1^-$ $(n-1)$ times, we denote this saddle connection by $\eta_3$. It follows that $\Sig$ is decomposed into three vertical cylinders, which will be denoted by $C_i^\vee, i=1,2,3,$ where

\begin{itemize}

\item[.] $C_1^\vee$ is bounded by $\eta_i^\pm$,

\item[.] $C_2^\vee$ is bounded by $\eta_1^+\cup\eta_2^+$ and $\eta_1^-\cup \eta_2^-$,

\item[.] $C_3^\vee$ is bounded by $\eta_2^+\cup\eta_3$ and $\eta_2^-\cup \eta_3$ ($\eta_3$ bounds $C_3^\vee$ from both sides).

\end{itemize}

\noindent Fix $\ell_2=h_2=1$, and let $m_i$ denote the modulus of $C_i$, and $m^\vee_i$ denote the modulus of $C^\vee_i, i=1,2,3$. Set

\begin{itemize}
\item[.]$\DS{a=\frac{m_1}{m_2}=\frac{h_1}{\ell_1}}$,

\item[.] $\DS{b=\frac{m_3}{m_2}=\frac{h_3}{\ell_3}}$,

\item[.] $\DS{c=\frac{m^\vee_2}{m^\vee_1}=\frac{h_3(1-\ell_1)}{(h_3+1)(\ell_1+\ell_3-1)}}$,

\item[.] $\DS{d=\frac{m^\vee_2}{m^\vee_3}=\frac{n^2(h_1+1)(1-\ell_1)}{(h_3+1)\ell_1}}$.

\end{itemize}

\noindent Let $\eta$ denote the saddle connection in $C_1$ which corresponds to the vector $\DS{(-\frac{\ell_1}{n}, h_1)}$. We have a pair of homologous saddle connections $\sigma_1^\pm$ in $T^\vee_3$ such that $\DS{\sigma_1^+=\sigma_1^-=\delta_1^-+\eta_2^-+\eta}$ in $H_1(\Sig, \Z)$.  Note that $\sigma^\pm_1$ bound a cylinder containing $\eta$. Similarly, we have in $T^\vee_2$ a pair of homologous saddle connections $\sigma_2^\pm$ such that $\DS{\sigma_2^+=\sigma_2^-=\delta_2^++\eta_1^++\eta_2^+}$. We have

$$V(\sigma_1^+)=(\frac{n-1}{n}\ell_1, h_1+1), \;  V(\sigma_2^+)=(1-\ell_1, h_3+1).$$

\noindent \underline{\bf Claim 1:} \textit{$V(\sigma_1^+)$ and $V(\sigma_2^+)$ are collinear if and only if $d=n(n-1)$.}\\

\noindent \textbf{Proof:} The fact that  $ V(\sigma_1^+)$ is parallel to $V(\sigma_2^+)$ is equivalent to

\begin{eqnarray*}   
  \frac{n-1}{n}\frac{\ell_1}{1-\ell_1} & = & \frac{1+h_1}{1+h_3}\\   
  \Longleftrightarrow  \frac{1-\ell_1}{1+h_3}\frac{1+h_1}{\ell_1} & = & \frac{n-1}{n}  \\  
   \Longleftrightarrow   \frac{d}{n^2}& = & \frac{n-1}{n} \\   
   \Longleftrightarrow  d & =  & n(n-1).
\end{eqnarray*}

\carre

Clearly, the surface $\Sig$ is completely determined by the values of $(h_1,\ell_1,h_3,\ell_3)$. We will find some values of $(h_1,\ell_1,h_3,\ell_3)$ such that the vertical and horizontal directions are parabolic, \ie $a,b,c,d$ are rational numbers, and $\Sig$ admits a decomposition into three cylinders in the direction $V(\sigma_1^\pm)$ whose moduli are independent over $\Q$. For this purpose, we fix $n\in \N$, $a,b,c$ in $\Q$, and $d=n(n-1)$, then we compute $(h_1,\ell_1,h_3,\ell_3)$ as functions of $(n,a,b,c)$.\\

\noindent First, observe that since $\DS{c=\frac{m_2^\vee}{m_1^\vee}}$ is a rational number, the vector $V(\sigma_2^+)$ is parallel to a vector in the lattice associated to $T^\vee_1$. It follows that $\Sig$ is decomposed into three cylinders in the direction $V(\sigma_1^+)$. We denote these cylinders by $C'_i, i=1,2,3$, where

\begin{itemize} 

\item[.] $C'_1$ is the cylinder containing $\eta$, and bounded by $\sigma_1^\pm$,

\item[.] $C'_2$ is the cylinder containing $\eta_2^\pm$ , and bounded by $\sigma_1^+\cup\sigma_2^+$ and $\sigma_1^-\cup\sigma_2^-$,

\item[.] $C'_3$ is the complement of $C'_1\cup C'_2$, which is bounded by $\sigma_2^+\cup \sigma_3$ and $\sigma_2^-\cup\sigma_3$, where $\sigma_3$ is a saddle connection parallel to $\sigma_2^\pm$, and bounds $C'_3$ from both sides.

\end{itemize}

\noindent Let $m'_i, i=1,2,3,$ denote the modulus of $C'_i$. \\

\noindent \underline{\bf Claim 2:} \textit{Suppose that $c\in \N$, then we have}

\begin{eqnarray*}
\frac{m'_1}{m'_2} & = & \frac{1}{(n-1)^2}(na\ell_1+1)(\frac{n}{\ell_1}-1),\\
\frac{m'_3}{m'_2} & = & \frac{b}{nc^2}(n-\ell_1).
\end{eqnarray*}

\noindent \textbf{Proof:} Let $h'_i, \ell'_i, i=1,2,3$, denote the height and the width of $C'_i$. We have

\begin{equation*}
\frac{h'_1}{h'_2}=\frac{\Aa(C'_1)}{\Aa(T^\vee_3)-\Aa(C'_1)}=\frac{(h_1+1)\ell_1-((n-1)/n)\ell_1}{((n-1)/n)\ell_1}=\frac{1+nh_1}{n-1}.
\end{equation*}

\noindent Using the vertical projection onto the horizontal axis, we have

\begin{equation*}
\frac{\ell'_2}{\ell'_1}= 1+\frac{1-\ell_1}{((n-1)/n)\ell_1}=\frac{n-\ell_1}{(n-1)\ell_1}.
\end{equation*}

\noindent Therefore 

$$\frac{m'_1}{m'_2}=\frac{h'_1}{h'_2}\frac{\ell'_2}{\ell'_1}=\frac{1+nh_1}{(n-1)}\frac{n-\ell_1}{(n-1)\ell_1}=\frac{1}{(n-1)^2}(1+na\ell_1) \frac{n-\ell_1}{\ell_1}.$$

\noindent Rescaling so that $V(\delta_3)=(1,0)$ and $V(\eta_1)=(0,1)$, that is $m^\vee_1=1$. Since $\DS{\frac{m^\vee_2}{m^\vee_1}=c}$, we deduce that $V(\sigma_2^+)$ is collinear with the vector $(c,1)$. When $c$ is an integer, in the standard torus $\R^2/\Z^2$, the number of intersection points of the simple closed geodesics corresponding to the vectors $(c,1)$ and $(0,1)$ is given by $|(c,1)\wedge (0,1)| =c$. Therefore, the saddle connection $\sigma_3$ crosses $\eta_1^+$ $c$ times. Using the projection along $V(\sigma_1^\pm)$ onto the vertical axis we have $\DS{ \frac{h'_3}{h'_2}=\frac{h_3}{c} }$. Then using the vertical projection onto the horizontal axis, we have $\DS{ \frac{\ell'_2}{\ell'_3}=\frac{1-\ell_1/n}{c\ell_3}}$. Therefore,

\begin{equation*}
\frac{m'_3}{m'_2}=\frac{h'_3}{h'_2}\frac{\ell'_2}{\ell'_3}=\frac{1}{nc^2}\frac{h_3}{\ell_3}(n-\ell_1)=\frac{b}{nc^2}(n-\ell_1).
\end{equation*}

\carre

Given $n\in\N, n>1$, $a,b \in \Q, a>0,b>0$, and $c \in \N, c>0$, set

$$P(X)=\frac{nc}{n-1}(X+a)(\frac{n}{(n-1)b}(X-1)(X+a)-\frac{X}{b}-X+1)-\frac{n}{n-1}(X-1)(X+a)+X.$$

\noindent \underline{\bf Claim 3:} \textit{Suppose that $P(X)$ is irreducible over $\Q$, and has a real root $\alpha$ satisfying }

\begin{itemize}

\item[.] $\alpha>1$,

\item[.] $\DS{\frac{n(\alpha-1)(\alpha+a)}{(n-1)\alpha}>1}.$

\end{itemize}

\noindent \textit{Then by taking $\DS{\ell_1 =   \frac{1}{\alpha}, h_1  =  \frac{a}{\alpha},\ell_3  =  \frac{1}{b}(\frac{n(\alpha-1)(\alpha+a)}{(n-1)\alpha} -1), h_3  =  \frac{n(\alpha-1)(\alpha+a)}{(n-1)\alpha} -1}$, the construction above gives us a Thurston-Veech surface with trace field of degree $3$ over $\Q$ for which the moduli of the three cylinders in the direction $V(\sigma_1^\pm)$ are independent over $\Q$. Consequently, this surface is generic in $\H^\mathrm{hyp}(4)$ by Corollary \ref{CorB}}.\\

\noindent \textbf{Proof:} From the choice of $\ell_1, h_1,\ell_3, h_3$, we only need to check that 

\begin{itemize}
\item[a)] $\DS{\ell_1+\ell_3-1>0}$

\item[b)] $\DS{\frac{m^\vee_2}{m^\vee_1}=\frac{1-\ell_1}{1+h_3}\frac{h_3}{\ell_1+\ell_3-1}=c}$

\item[c)] $\DS{\frac{m^\vee_2}{m_3^\vee}=\frac{n^2(1-\ell_1)(1+h_1)}{\ell_1(1+h_3)}=n(n-1)}$.

\end{itemize}

\noindent Condition a) is satisfied since we have

\begin{eqnarray*}
\ell_1+\ell_3-1 & = & \frac{1}{\alpha}+\frac{n}{b(n-1)}\frac{(\alpha-1)(\alpha+a)}{\alpha}-\frac{1}{b}-1\\
                & = & \frac{n-1}{cn(\alpha+a)}(\frac{n}{n-1}\frac{(\alpha-1)(\alpha+a)}{\alpha}-1) >0.
\end{eqnarray*}

\noindent Conditions b) and c) follow immediately from the fact that $\alpha$ is a root of $P$. To see that the trace field of $\Sig$ is of degree $3$, remark that we have $m_i \in \Q, i=1,2,3,$ but 

$$m_3^\vee=\frac{\ell_1}{n^2(1+h_1)}=\frac{1}{n^2(\alpha+a)}.$$

\noindent is an algebraic number of degree $3$ over $\Q$. From Claim 2, we have

\begin{eqnarray*}
\frac{m'_1}{m'_2} & = & \frac{1}{(n-1)^2}(1+na\ell_1)(\frac{n}{\ell_1}-1)\\
 & = & \frac{1}{(n-1)^2}(1+\frac{na}{\alpha})(n\alpha-1)\\
 & = & \frac{n}{(n-1)^2}\alpha + \frac{n^2a-1}{(n-1)^2} - \frac{na}{(n-1)^2}\frac{1}{\alpha}.
\end{eqnarray*}

\noindent and

\begin{eqnarray*}
\frac{m'_3}{m'_2} & = & \frac{b}{nc^2}(n-\ell_1)\\
& = & \frac{b}{c^2} - \frac{b}{nc^2}\frac{1}{\alpha}.
\end{eqnarray*}

\noindent Since $\DS{\alpha, 1, \frac{1}{\alpha}}$ are independent over $\Q$, it follows that $m'_1,m'_2,m'_3$ are independent over $\Q$. \carre

\subsection{Numerical examples}

Here below are some explicit examples of Thurston-Veech surfaces obtained from the construction above which satisfy the condition of Corollary \ref{CorB}. Here $\tilde{P}(X)$ is a polynomial proportional to $P(X)$ with coefficients in $\Z$.

\begin{center}
\begin{tabular}{|c|c|c|c|c|}
  \hline
  $(n,a,b,c)$ & $\tilde{P}(X)$ & $\alpha$ & $\ell_1$ & $\ell_3$ \\
  \hline
  $(4,1,10,5)$ &  $8X^3-70X^2-5X+64 $& $\approx 8.716407 $ & $\approx 0.114726 $ & $\approx 1.046891$ \\
  $(5,2,10,3)$ & $15X^3-127X^2-152X+260$  & $\approx 9.352026$ & $\approx 0.106929$ & $\approx 1.167271$ \\
  $(5,1/5,2,1)$ & $25X^3-115X^2+83X+15$ & $\approx 3.643625$ & $\approx 0.274452$ & $ \approx 1.242959$ \\
  $(5,1/2,5,1)$ & $20X^3-176X^2+121X+75$ & $\approx 7.983332$ & $\approx 0.125261$ & $\approx 1.655175$\\
  $(2,1,6,1)$ & $2X^3-11X^2+10$ & $\approx 5.323574$ & $\approx 0.187844$ & $\approx 1.545243$\\
  $(2,2,9,2)$ & $8X^3-34X^2-53X+76$ & $\approx 5.175414$ & $\approx 0.193221$ & $\approx 1.175327$\\
  \hline
\end{tabular}

\end{center}

\end{document}